\documentclass[a4paper,11pt]{amsart}

\usepackage{amsmath, amssymb, amsthm, fullpage}
\usepackage{url}

\usepackage{graphicx}

\usepackage{breqn} 

\usepackage[
margin=0.04\textwidth
]{caption} 

\DeclareMathOperator{\Rat}{Rat}
\DeclareMathOperator{\Disc}{Disc}
\DeclareMathOperator{\Res}{Res}
\DeclareMathOperator{\Int}{int}
\DeclareMathOperator{\PSL}{PSL}
\DeclareMathOperator{\Spec}{Spec}
\DeclareMathOperator{\Proj}{Proj}
\DeclareMathOperator{\Fix}{Fix}
\DeclareMathOperator{\SL}{SL}
\DeclareMathOperator{\GL}{GL}
\newcommand{\Id}{\mathrm{Id}}
\newcommand{\bA}{\mathbb{A}}
\newcommand{\bP}{\mathbb{P}}
\newcommand{\bQ}{\mathbb{Q}}
\newcommand{\bR}{\mathbb{R}}
\newcommand{\bC}{\mathbb{C}}
\newcommand{\bL}{\mathbb{L}}
\newcommand{\bN}{\mathbb{N}}
\newcommand{\bZ}{\mathbb{Z}}
\newcommand{\rM}{\mathrm{M}}
\newcommand{\cB}{\mathcal{B}}
\newcommand{\cC}{\mathcal{C}}
\newcommand{\cS}{\mathcal{S}}
\newcommand{\cH}{\mathcal{H}}
\newcommand{\cO}{\mathcal{O}}
\newcommand{\an}{\mathrm{an}}

\numberwithin{equation}{section}
\theoremstyle{plain}
\newtheorem{theorem}{Theorem}[section]

\newtheorem{proposition}[theorem]{Proposition}
\newtheorem{mainth}{Theorem}
\newtheorem{maincoro}{Corollary}

\theoremstyle{definition}

\newtheorem{notation}[theorem]{Notation}

\newtheorem*{acknowledgement}{Acknowledgement}
\theoremstyle{remark}
\newtheorem{remark}[theorem]{Remark}

\begin{document}
\title[Parabolic bifurcation loci in the spaces of rational functions]{
Parabolic bifurcation loci in the spaces of rational functions}
\author[Y\^usuke Okuyama]{Y\^usuke Okuyama
}
\address{
Division of Mathematics,
Kyoto Institute of Technology, Sakyo-ku,
Kyoto 606-8585 Japan}
\email{okuyama@kit.ac.jp}

\date{\today}

\subjclass[2010]{Primary 37F45; Secondary 32H50}
\keywords{parabolic bifurcation locus, parabolic bifurcation hypersurface, 
dynatomic polynomial, multiplier polynomial, cyclotomic polynomial,
complex dynamics, non-archimedean dynamics, potential theory}

\begin{abstract}
We give a geometric description of the parabolic bifurcation locus
in the space $\Rat_d$ of all rational functions on $\bP^1$ 
of degree $d>1$, generalizing the study by Morton and Vivaldi 
in the case of monic polynomials. The results are new even for quadratic
rational functions.
\end{abstract}

\maketitle
 
\renewcommand{\thefootnote}{\fnsymbol{footnote}}
\footnote[0]{This is the Accepted Manuscript version of an article accepted for publication in Nonlinearity. IOP Publishing Ltd is not responsible for any errors or omissions in this version of the manuscript or any version derived from it. The Version of Record is available online at 
\url{https://doi.org/10.1088/1361-6544/ac8041}.}
\renewcommand{\thefootnote}{\arabic{footnote}}

\section{Introduction}\label{sec:intro}

In studying a (topologically parametrized)
family of (discrete) dynamical systems on a (common) phase space, 
once a stability notion for this family is fixed, 
one of main topics is classifying
the (in)stability phenomena 
under a perturbation/variation of the dynamical system.
We also call an instability phenomenon a bifurcation.
For example, we say a family of dynamical systems
is structurally stable at a dynamical system $f$
if any small perturbation of $f$ is conjugate to the initial $f$ 
by some automorphism of the phase space. The structural stability
is an open condition, but the structurally stable locus in the parameter space
is not necessarily dense even if it is non-empty. 
An obvious obstacle to the structural stability 
(or a possible unstable limit of structurally stable dynamics)
is (informally) a collision between two distinct periodic points 
(possibly belonging to a same cycle) of a dynamical system.
In other (but still naive) words, by perturbing 
such an unstable dynamical system appropriately, 
some multiple cycle of it splits into more than one cycles (having
the same period) or 
turns to a single but simple cycle (having a greater period), and we call such a (classical) bifurcation phenomenon a parabolic bifurcation.

Even for the (iterations of the monic centered 
complex) quadratic polynomial family 
\begin{gather*}
 z\mapsto z^2+t\quad\text{on }\bC 
\end{gather*}
(holomorphically) parametrized by $t\in\bC$ (which is a toy model of 
the dynamical moduli space $\rM_d$ of degree $d>1$ rational functions on $\bP^1$ below), the theory of stability/bifurcation 
is much richer than what would be expected. 
For each individual 
$P_t(z)=z^2+t$,
which is regarded as an endomorphism of $\bC$ (or 
the projective line $\bP^1(\bC)=\bC\cup\{\infty\}$
equipped with the chordal metric),
the phase space $\bC$ 
is partitioned into a non-equicontinuity part
(the chaotic locus in the sense of Devaney) 
and an equicontinuity part (the 
region of normality in the sense of Montel) for the sequence of
the iterations $P_t^n=P_t^{\circ n}=P_t\circ\cdots\circ P_t$ 
($n$ times), $n\in\bN$, which are respectively called the Julia set $J_t$
and the Fatou set $F_t$ of $P_t$ and 
are both totally invariant under $P_t$ (for complex dynamics,
we refer to Milnor \cite{Milnor3rd}).
We say the family $(P_t)_{t\in\bC}$ 
is $J$-stable at the parameter $t=t_0\in\bC$ 
if the (compact set valued) function $t\mapsto J_t$ 
is continuous at $t=t_0$ in the Hausdorff topology.
By a general theory due to Ma\~n\'e--Sad--Sullivan, 
Lyubich, and McMullen (\cite{MSS,Lyubich83stability},
and see also \cite[\S4]{McMullen:renorm})
on a holomorphic family of 
rational functions on $\bP^1(\bC)$ of (common) degree $d>1$, 
the $J$-unstable (or bifurcation) locus 
$\cB_2$ in the parameter space $\bC$ coincides 
with the (fractal-shaped) boundary $\partial\cC_2$ of the so called 
Mandelbrot set 
\begin{gather*}
\cC_2=\Bigl\{t\in\bC:\limsup_{n\to\infty}|P_t^n(0)|<+\infty\Bigr\} 
\end{gather*}
(where we note that
the point $z=0$ is the unique critical point of $P_t$ in the phase space $\bC$),
the $J$-stable locus $\cS_2=\bC\setminus\partial\cC_2$ 
in the parameter space $\bC$ is (open and) dense, and
the difference between $\cS_2$ and the structurally stable locus
in the parameter space $\bC$ of $(P_t)_{t\in\bC}$ is the set of all the
(so called) superattracting parameters $t$, that is, 
the roots $t$ in $\bC$ of the polynomial
\begin{gather*}
 P_t^n(0)=\bigl(P_t^n(z)-z\bigr)\bigl|_{z=0}\,\in\bZ[t] 
\end{gather*} 
for some $n\in\bN$. We say a parameter $t\in\bC$ is hyperbolic if the (restricted dynamics) 
$P_t:J_t\to J_t$ is uniformly expanding. The superattracting parameters are hyperbolic ones.
The hyperbolicity locus $\cH_2$ is clopen in $\cS_2$. Conjecturally, $\cH_2=\cS_2$.

The bifurcation locus $\cB_2=\partial\cC_2$ also densely contains 
various kinds of distinguished parameters $t$ (at least nearly)
algebraically defined. Most importantly,
$\cB_2$ densely contains all the parabolic bifurcation
parameters mentioned at the beginning;
for example, at the parameter $t=1/4$, $P_{1/4}(z)-z=z^2-z+1/4$ 
has the multiple root $z=1/2$.
Formally, the parabolic bifurcation (collision)
parameters $t$ for the quadratic polynomial family $(P_t)_{t\in\bC}$
should have been defined by the roots $t$ in $\bC$ of the discriminants 
\begin{gather*}
 \Disc\bigl(P_t^n(z)-z\bigr)\in\bZ[t]
\end{gather*}
of $P_t^n(z)-z\in(\bZ[t])[z]$, $n\in\bN$, which contains all the roots $t$ in $\bC$ of the resultant
\begin{gather*}
 \Res\bigl(P_t^n(z)-z,P_t^\ell(z)-z\bigr)\in\bZ[t]
\end{gather*}
between $P_t^n(z)-z$ and $P_t^\ell(z)-z\in(\bZ[t])[z]$, 
$n\in\bN$ and $\ell\in\bN$ dividing $n$ and less than $n$.
Another example of an (almost) algebraic parameter 
$t\in\cB_2$ is a Misiurewicz one, for which $z=0$ is 
in $J_t$ as well as is (strictly) preperiodic, 
i.e.\ $t$ is a root of 
$P_t^n(0)-P_t^m(0)\in\bZ[t]$ 
for some integers $n>m>0$, and it is observed by a computer drawing
(and shown by renormalization) that $\cB_2$ near a Misiurewicz
parameter $t$ is asymptotically similar to $J_t$ 
near the (so called postcritical) orbit $z=P_t^n(0)$ for $n\gg 1$ 
(since \cite{Tan90}; for a recent potential
theoretic development, see \cite{AGMV19}). 

Although both parabolic bifurcation parameters
and Misiurewicz parameters are dense in $\cB_2$, 
it has not 
been completely understood how they 
are equidistributed towards the harmonic measure on $\cB_2$ quantitatively 
(in terms of electrostatistic or potential theory,
since \cite{Levin90}; for recent developments, 
see e.g.\ \cite{BassanelliBerteloot11,BG14,Gauthier13,GOV16}). 
Other interesting species also live in $\cB_2$ (see e.g.\ \cite{Widom83,McSiegel,DLN20}).

\begin{figure}[h]
\includegraphics[bb=0 0 999 999,width=0.35\textwidth]{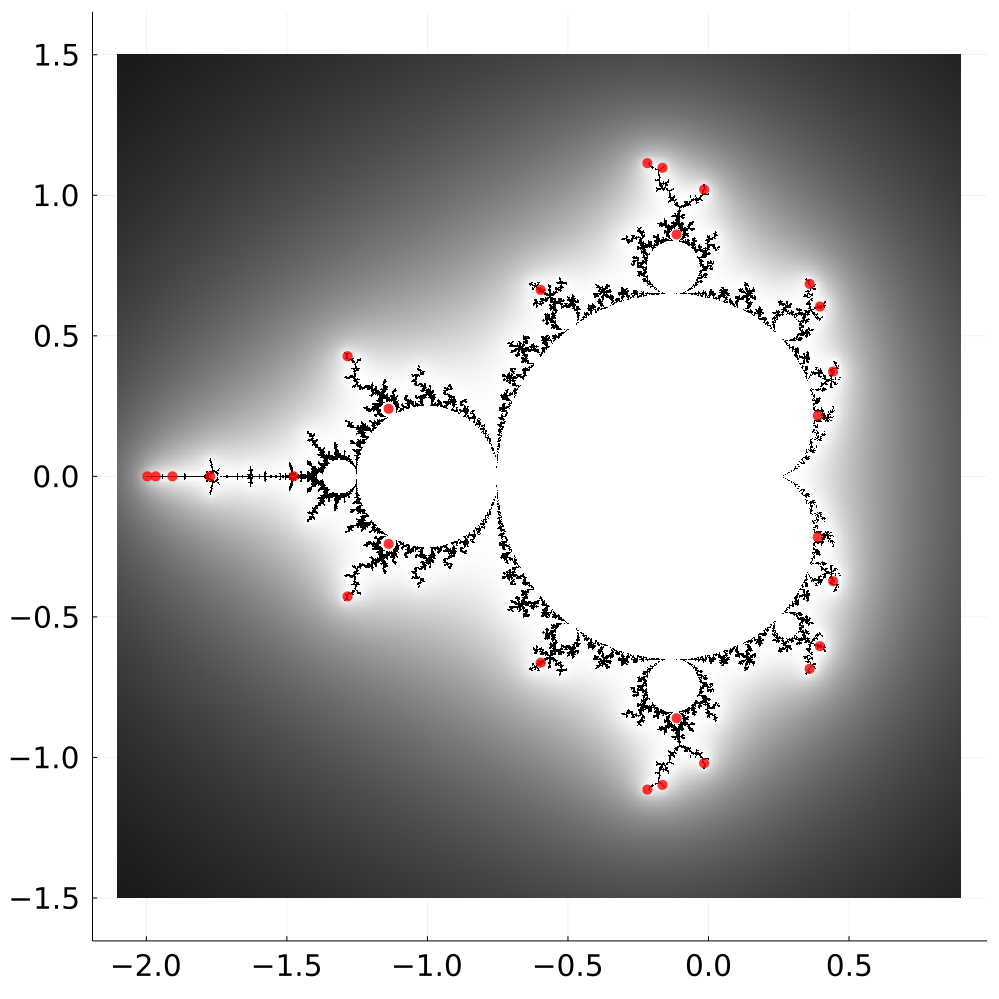}
\includegraphics[bb=0 0 999 999,width=0.35\textwidth]{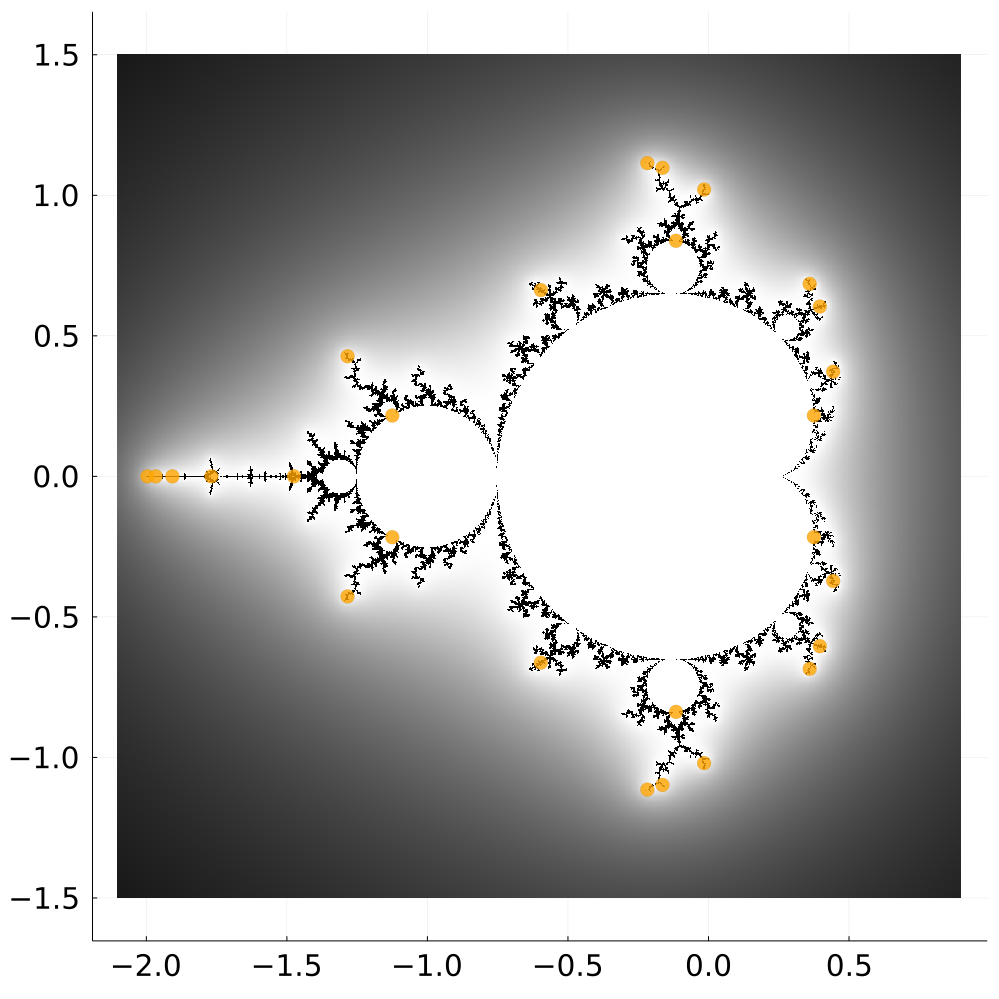}
\caption{The Mandelbrot set $\cC_2$ (white),
the superattracting parameters $t$ (red dots, living in $\Int\cC_2\subset\cS_2$),
and the parabolic bifurcation 
parameters $t$ (orange dots, 
living in $\cB_2=\partial\cC_2$) for the (formally exact) periods $n=6$; 
those two computer pictures exhibit that
the distribution of the superattracting parameters
and that of the parabolic bifurcation parameters are similar,
and indeed the averaged mass distributions/counting measures
of both converge to the (same) harmonic measure/equilibrium mass distribution 
$\mu_{\cB_2,\infty}$ on $\cB_2$ with pole $\infty$
as the (formally exact) period $n\to\infty$. 
However, no (non-trivial) speed estimates of convergence of the parabolic bifurcation parameters towards $\mu_{\cB_2,\infty}$ have been known, while that of the superattracting ones
has been known quite precisely (see e.g.\ \cite{GOV16}); in other 
words, no non-trivial quantitative asymptotic estimates of the difference 
between the (averaged distributions of the) superattracting and 
parabolic bifurcation parameters as the period $n\to\infty$
have been known.}
\end{figure}

The study of arithmetic dynamics
(dynamics of rational functions on $\bP^1$ defined over 
a number field e.g.\ $\bQ$ which has
the $p$-adic (non-archimedean) norm $|\cdot|_p$ 
for each prime $p$ as well as the Euclidean (archimedean) norm $|\cdot|_\infty$)
plays a more and more important role
in the study of complex dynamics (i.e.\ dynamics of
rational functions on $\bP^1$ defined over $\bC$); 
it might be surprising that
the theory of heights in arithmetic is closely related to
electrostatistic or potential theory (since Baker--Rumely, Chambert-Loir, and Favre--Rivera-Letelier \cite{BR10,ChambertLoir06,FR06}).
The ($J$-)stability and bifurcation in non-archimedean setting
are more subtle (e.g.\ see \cite{SilvermanLambda}
for a non-archimedean counterpart to the
Ma\~n\'e--Sad--Sullivan and Lyubich theory).
Morton and Vivaldi \cite{MortonVivaldi95} studied,
towards a general bifurcation theory
in (merely)
algebraic dynamics,
the parabolic bifurcation locus (or hypersurfaces)
in the ($d$-dimensional) parameter space $R^d$ of monic polynomials
\begin{gather*}
 z^d+b_{d-1}z^{d-1}+\cdots+b_1z+b_0,\quad
(b_0,\ldots,b_{d-1})\in R^d
\end{gather*}
of degree $d$, $d>1$, 
defined over an integral domain $R$
(e.g.\ the ring $\bZ$).
Our aim is to contribute to the study of stability and bifurcation
in algebraic(, complex, and non-archimedean) dynamics
by developing the Morton--Vivaldi-type theory in 
the space $\Rat_d$ of degree $d$ rational functions on $\bP^1$
(and also in the dynamical moduli space $\rM_d$ of them), so in particular
by giving a geometric description of
the parabolic bifurcation hypersurfaces in $\Rat_d$ and $\rM_d$.

\subsection*{Organization of the paper}
In Section \ref{sec:mainth}, we introduce the necessary notions
and state the main results (Theorem \ref{th:exact} and Corollary \ref{th:geometric}).
In Section \ref{sec:background},
we recall some complex dynamical and potential theoretic calculus on 
$\bP^1=\bP^1(\bC)$, and
also some details on the M\"obius function/inversion and
on homogeneous resultants/discriminants.
In Section \ref{sec:proofself}, we show Theorem \ref{th:exact}.
In Section \ref{sec:quadrat}, we conclude with a few computations in the case
where $d=2$.

\section{Main results}\label{sec:mainth}

\subsection{The space $\Rat_d$ of degree $d$ rational functions on $\bP^1$}\label{sec:scheme}
The space of rational functions on $\bP^1$ is by no means
(a little) more complicate than that of polynomials.
For each integer $d\ge 0$, over the complex number field $\bC$, 
the set $\Rat_d(\bC)$ of all rational functions
\begin{gather*}
f(z)=P(z)/Q(z),\quad P,Q\in\bC[z]\text{ are coprime, and }\max\{\deg P,\deg Q\}=d 
\end{gather*}
on $\bP^1(\bC)$ of degree $d$ is regarded as a complex manifold
$\bP^{2d+1}(\bC)\setminus V_d(\bC)$,
by identifying $f$ with the ratio 
\begin{gather*}
 [a:b]=[a_0:\cdots:a_d:b_0:\cdots:b_d]\in\bP^{2d+1}(\bC),\quad a=(a_0,\ldots,a_d),\,b=(b_0,\ldots,b_d)\text{ for short},
\end{gather*}
of the $2d+2$ coefficients $a_0,\cdots,a_d,b_0,\ldots,b_d\in\bC$ of the denominator and the numerator
\begin{gather*}
 Q(z)=\sum_{j=0}^da_jz^j\quad\text{and}\quad P(z)=\sum_{k=0}^db_kz^k 
\end{gather*}
of $f$, but by removing a (complex analytic) 
hypersurface $V_d(\bC)$ from 
the ambient $2d+1$ dimensional projective space
$\bP^{2d+1}(\bC)$ so that $P,Q$ are coprime.
We note that $\Rat_0(\bC)=\bP^1(\bC)$
(and $V_0(\bC)=\emptyset$), and that 
the $\Rat_1(\bC)$, the M\"obius transformation/projective coordinate 
change group on $\bP^1(\bC)$, is $\cong\PSL(2,\bC)$. 

From now on, we assume $d>1$.

\begin{notation}
 For a (commutative) ring $R$, $R^*$ denotes the unit (or invertible element) group.
\end{notation}
\begin{notation}
  For each ring $R$ and each $d_0\in\bN\cup\{0\}$, 
 the set 
 of all homogeneous polynomials in $R[X_0,\ldots,X_n]$
 of degree $d_0$ is denoted by $R[X_0,\ldots,X_n]_{d_0}$. 
 Here,
 $H\in R[X_0,\ldots,X_n]$ is said to be homogeneous if
 $H$ is the sum of monomials in $X_0,\ldots,X_n$
 the sum of the powers of $X_0,\ldots,X_n$ of any of which
 (is identical and) equals $\deg H$.
\end{notation}

Then
we (abstractly) introduce the space
of degree $d$ rational functions on $\bP^1$ 
\begin{gather*}
 \Rat_d:=\bP^{2d+1}_{\bZ}\setminus V_d
\end{gather*}
as a scheme over the scheme $\Spec\bZ$, where 
$V_d:=V(\rho_d)$ is the hypersurface in $\bP^{2d+1}_{\bZ}$
defined by (the zeros in $\bP^{2d+1}_{\bZ}$ of)
the $d$-th homogeneous 
``resultant'' 
form
\begin{gather*}
\rho_d(a,b):=\det\begin{pmatrix}
	    a_0 & \cdots & a_{d-1} & a_d  &        &    \\
	    \text{\huge{0}}    & \ddots & \vdots  & \vdots & \ddots & \text{\huge{0}}   \\
	        &        & a_0     & a_1  & \cdots & a_d\\
	    b_0 & \cdots & b_{d-1} & b_d  &        &    \\
	    \text{\huge{0}}  & \ddots & \vdots  & \vdots & \ddots &  \text{\huge{0}}  \\
	        &        & b_0     & b_1  & \cdots & b_d\\ 
	   \end{pmatrix}\in\bZ[a_0,\ldots,a_d,b_0,\ldots,b_d]_{2d}\label{eq:resultant}
\end{gather*}
on $\bP^{2d+1}_{\bZ}$ (or on $\bA^{2d+2}_{\bZ}$), 
still writing $a=(a_0,\ldots,a_d)$ and $b=(b_0,\ldots,b_d)$ for short.

We introduced $\Rat_d$ as a (not only quasi-projective but also affine)
scheme (defined over the ring $\bZ$),
which could be thought of as some topological space similar to an algebraic variety 
(defined over an algebraically closed field, e.g., $\bC$)
but still equipped with the ring 
\begin{gather}
\bZ(\Rat_d)=\Bigl\{\rho_d^{-\frac{\deg P}{2d}}\cdot P:P\in\bZ[a,b]\setminus\{0\}
\text{ is homogeneous and }(2d)|(\deg P)\Bigr\}\cup\{0\}\label{eq:regular}
\end{gather}
of regular functions on $\Rat_d$; see \cite[\S 1]{Silverman98} or \cite[\S 4.3]{SilvermanDynamics}
(we refer to Hartshorne \cite[\S1-\S2]{Hartshorne77} for schemes/algebraic varieties).
Since the above $\rho_d(a,b)$ is irreducible in $\bZ[a,b]$
(see \cite[Chapter 13]{BGZ08}) and the ring $\bZ[a,b]$ is a 
uniquely factorization domain (UFD), we have
\begin{gather}
 (\bZ(\Rat_d))^*=\bZ^*(=\{\pm 1\})\label{eq:regularunit} 
\end{gather}
(see Proposition \ref{th:unit} below). 
On the other hand, we also denote by $\cO_{\Rat_d(\bC),\an}(D)$
the ring of complex analytic functions on each non-empty open subset $D$ 
in the complex manifold $\Rat_d(\bC)$.

\begin{remark}
 We would do calculus only on the complex manifold 
 $\Rat_d(\bC)$ to avoid too many technicalities
 while, for the formulation and the proof of our main result,
 the above scheme-theoretic introduction of $\Rat_d$ (and $\rM_d$ below) is 
 not only natural but also necessary.
\end{remark}

\subsection{Variation of periodic points/cycles}\label{sec:dynatomic}
Computation of the iterations of rational functions on $\bP^1$ is also by no
means more complicate than that of polynomials.
As in the previous section, we identify each 
$[a:b]=[a_0:\cdots:a_d:b_0:\cdots:b_d]\in\bP^{2d+1}_{\bZ}\setminus V_d$ with a rational function
\begin{gather*}
f(z)=f_{[a:b]}(z)=f_{[a_0:\cdots:a_d:b_0:\cdots:b_d]}(z)
=\frac{\sum_{j=0}^db_jz^j}{\sum_{i=0}^da_iz^i}
\quad\text{on }\bP^1_{\bZ}(=\Proj(\bZ[X,Y])), 
\end{gather*}
where $z=Y/X$ is an affine coordinate of $\bP^1_{\bZ}$,
and then the pair of the 
homogenizations 
\begin{multline}
 F(X,Y)=F(X,Y;a,b)=\bigl(F_0(X,Y;a,b),F_1(X,Y;a,b)\bigr)\\
=\bigl(F_0(X,Y),F_1(X,Y)\bigr)=\Bigl(\sum_{i=0}^da_iX^{d-i}Y^i,\sum_{j=0}^db_jX^{d-j}Y^j\Bigr)\in\bigl((\bZ[a,b]_1)[X,Y]_d\bigr)^2\label{eq:lift}
\end{multline}
of the denominator and the numerator of $f=f_{[a:b]}$(, which is
identified with the point $(a,b)\in\bA^{2d+2}_{\bZ}$, and is unique up to the
``multiplication by a non-zero constant'' action on $\bA^{2d+2}_{\bZ}$
of the multiplication group scheme $\mathbb{G}_m$,) is an endomorphism of the
affine plane $\bA^2_{\bZ}(=\Spec(\bZ[X,Y]))$ of degree $d$ and is
called a lift (to $\bA^2_{\bZ}$) of $f=f_{[a:b]}$.  Such a lifting $f$
(to $\bA^2_{\bZ}$) is useful in computation, noting that for each
$n\in\bN$, the $n$-th iteration $F^n=F\circ\cdots\circ F$ ($n$ times) of $F$, which we write as
\begin{gather}
 F^n=\bigl(F_0^{(n)}(X,Y;a,b),F_1^{(n)}(X,Y;a,b)\bigr)\in\bigl((\bZ[a,b]_{(d^n-1)/(d-1)})[X,Y]_{d^n}\bigr)^2,\label{eq:iterlift} 
\end{gather}
is a lift (to $\bA^2_{\bZ}$)
of the $n$-th iteration $f^n=f\circ\cdots\circ f$ ($n$ times) in $\Rat_{d^n}$ of $f=f_{[a:b]}$. 

Letting $\Omega$ be an algebraically closed field of characteristic $0$, e.g.\ $\bC$, 
for each individual rational function
$f=f_{[a:b]}\in\Rat_d(\Omega)$ on $\bP^1(\Omega)$ of degree $d$
defined over $\Omega$,
an endomorphism $F=F(X,Y;a,b)\in(\Omega[X,Y]_d)^2$ of $\bA^2(\Omega)$,
where $(a,b)\in\Omega^{2d+2}$,
is also called a lift (to $\Omega^2$) of this $f=f_{[a:b]}$ and is unique up to a 
multiplication in $\Omega^*=\Omega\setminus\{0\}$.

Recall that we are interested in when $f$ has a
multiple cycle, that is, when or where a collision of distinct 
periodic points in $\bP^1$ of a rational function $f_{[a:b]}$ 
happens as the point $[a:b]$ varies in $\Rat_d$. 
The following notion from arithmetic is useful for our purpose.
For the M\"obius function $\mu:\bN\to\{0,\pm 1\}$ 
and the (additive/multiplicative) M\"obius inversion of a sequence 
in an abelian/commutative group, see Subsection \ref{th:mobius} below. 

\begin{notation}
For each $i\in\bN$, the $i$-th cyclotomic polynomial
\begin{gather*}
 C_i(T):=\prod_{k\mid i}(T^k-1)^{\mu(\frac{i}{k})}
=\prod_{\omega:\,\text{an $i$-th primitive root of unity (in }\overline{\bQ})}(T-\omega)\in\bZ[T]
\end{gather*}
in $T$ is the M\"obius inversion of the sequence $(T^i-1)_i$
(in $(\bQ(T))^*=\bQ(T)\setminus\{0\}$).
\end{notation}

First, let $\Omega$ be an algebraically closed field
of characteristic $0$.
For every individual $f\in\Rat_d(\Omega)$
and every $n\in\bN$, 
the set $\Fix(f^n)$ of all fixed points in $\bP^1(\Omega)$ of $f^n$,
i.e.\ all solutions of the equation 
$f^n(z)=z$
in $\bP^1(\Omega)$, is the zero locus in $\bP^1(\Omega)$ of the homogeneous polynomial
\begin{gather*}
 Y\cdot F_0^{(n)}(X,Y)-X\cdot F_1^{(n)}(X,Y)\in\Omega[X,Y]_{d^n+1}, 
\end{gather*}
where $F$ is a lift 
(to $\Omega^2$) of $f$
(recalling the notation in \eqref{eq:iterlift});
then we call
\begin{gather*}
 \Fix^*(f^n):=\Fix(f^n)\setminus\bigcup_{m:\,\mid n\text{ and }<n}\Fix(f^m) 
\end{gather*}
the set of all periodic points of $f$ in $\bP^1(\Omega)$ having 
the exact period $n$,
and denoting, for each $m\in\bN$ and each $z\in\Fix^*(f^m)$,
by $(f^m)'(z)(=(f^m)'(f(z))=\cdots=(f^m)'(f^{m-1}(z)))\in\Omega$ 
the multiplier of the cycle $z,f(z),\ldots,f^{m-1}(z)$ of $f$ in $\bP^1(\Omega)$
(even if $z=\infty$ for simplicity),
we call
\begin{gather}
\Fix^{**}(f^n):=\Fix^*(f^n)
\sqcup\bigsqcup_{m:\,\mid n\text{ and }<n}
\Bigl\{z\in\Fix^*(f^m):C_{\frac{n}{m}}\bigl((f^m)'(z)\bigr)=0\Bigr\}
\label{eq:decomp}
\end{gather}
the set of all periodic points of $f$ in $\bP^1(\Omega)$ having 
the {\em formally exact} period $n$ (which might be less studied but is more adequate 
for our purpose. 
The terminology ``formally exact period $n$'' is a slight modification of
the original ``formal period $n$'' in \cite[\S4.1]{SilvermanDynamics}). We note that 
\begin{gather*}
 \Fix^*(f^n)\subset\Fix^{**}(f^n)\subset\Fix(f^n),
\end{gather*}
the first inclusion in which is the equality for generic $f\in\Rat_d(\Omega)$.
\begin{remark}
 Our uses of the superscripts ``*'' and ``**'' are reverse to those in
 \cite[\S4.1]{SilvermanDynamics} but in this article, 
the less used ``**'' is attached to the less studied ``$\Fix^{**}(f^n)$''.
 The superscript ``**'' would be added to the notations coming from $\Fix^{**}(f^n)$.
\end{remark}

Next, for each $n\in\bN$, let $\Fix_n$ be
the hypersurface $V(\Phi_n)$ in the (fibered) product (scheme) 
\begin{gather*}
 \bP^1_{\Rat_d}:=\bP^1_{\bZ}\times\Rat_d
\end{gather*} 
defined by (the zeros in $\bP^1_{\Rat_d}$ of) the 
polynomial
\begin{multline*}
\Phi_n(X,Y;a,b)=\Phi_{F,n}(X,Y)\\
:=Y\cdot F_0^{(n)}(X,Y;a,b)-X\cdot F_1^{(n)}(X,Y;a,b)\in\bigl(\bZ[a,b]_{(d^n-1)/(d-1)}\bigr)[X,Y]_{d^n+1},
\end{multline*}
in terms of the notations in \eqref{eq:iterlift}.
This hypersurface $\Fix_n=V(\Phi_n)$ in $\bP^1_{\Rat_d}$
is regarded as a $d^n+1$-(multi)valued section of the projection
$\bP^1_{\Rat_d}\to\Rat_d$
and is encoding the variation of the $d^n+1$ periodic points in $\bP^1$
of $f=f_{[a:b]}$ having the exact periods dividing $n$ 
when the $[a:b]$ varies in $\Rat_d$, in working over an algebraically closed field
$\Omega$ of characteristic $0$.

Finally, for each $n\in\bN$, we also set the (additive) M\"obius inversion
\begin{gather*}
 d_n:=\sum_{m\mid n}\mu\Bigl(\frac{n}{m}\Bigr)(d^m+1)\quad(\text{so, e.g., }d_p=d^p-d\text{ for a prime number }p)
\end{gather*}
of the sequence $(d^n+1)_n$ (in the abelian group $\bZ$), and set
the M\"obius inversion
\begin{multline*}
\Phi_n^{**}(X,Y;a,b)=\Phi_{F,n}^{**}(X,Y)
:= \prod_{m\mid n}\bigl(\Phi_{F,m}(X,Y)\bigr)^{\mu(\frac{n}{m})}
\in
\begin{cases}
 (\bZ[a,b]_1)[X,Y]_{d+1} &\text{if }n=1,\\
 (\bZ[a,b]_{d_n/(d-1)})[X,Y]_{d_n} & \text{if }n>1,
\end{cases}
\end{multline*}
which is a priori a rational function but indeed a (homogeneous) polynomial, 
of the sequence $(\Phi_n)_n$
(in $((\bZ[a,b])(X,Y))^*=((\bZ[a,b])(X,Y))\setminus\{0\}$); 
letting $\Omega$ be an algebraically closed field
of characteristic $0$,
for each individual $f\in\Rat_d(\Omega)$, the homogeneous polynomial
$\Phi_{F,n}^{**}(X,Y)\in\Omega[X,Y]_{d_n}$ is called 
the $n$-th dynatomic polynomial for a lift $F$ (to $\Omega^2$) of $f$, and 
recalling \eqref{eq:decomp}, we indeed have
\begin{gather}
\Fix^{**}(f^n)=(\text{the zero locus in 
$\bP^1(\Omega)$ of $\Phi_{F,n}^{**}\in\Omega[X,Y]_{d_n}$}).\label{eq:formaldynam}
\end{gather}

For each $n\in\bN$, we have $\Phi_n=\prod_{m\mid n}\Phi_m^{**}$, and
letting $\Fix_n^{**}$ be the hypersurface $V(\Phi_n^{**})$ in $\bP^1_{\Rat_d}$
defined by (the zeros in $\bP^1_{\Rat_d}$ of) the polynomial $\Phi_n^{**}(X,Y;a,b)$, 
we also have 
\begin{gather*}
 \Fix_n=\bigcup_{m\mid n}\Fix_m^{**}\quad\text{in }\bP^1_{\Rat_d}.
\end{gather*}
The hypersurface
$\Fix_n^{**}=V(\Phi_n^{**})$ in $\bP^1_{\Rat_d}$ is regarded as
a $d_n$-(multi)valued section of the projection $\bP^1_{\Rat_d}\to\Rat_d$
and is indeed encoding the variation of the $d_n$ periodic points in $\bP^1$
of $f=f_{[a:b]}$ having the formally exact periods $n$ 
when the $[a:b]$ varies in $\Rat_d$, in working over an algebraically closed field
$\Omega$ of characteristic $0$.

See \cite[Theorem 4.4 and its Remark]{Silverman98} 
for a full account over $\bZ$, and also \cite[\S4.1]{SilvermanDynamics} 
over $\bQ$.

\subsection{Collision of periodic points: parabolic bifurcation locus in $\Rat_d$}\label{sec:main}
For the homogeneous resultant/discriminant of homogeneous polynomials
in $X,Y$, see Subsection \ref{sec:resultant}. 

First, for each $n\in\bN$, the homogeneous discriminant
$\Disc(\Phi_n^{**})$ of the homogeneous polynomial 
$\Phi_n^{**}\in(\bZ[a,b])[X,Y]_{d_n}$ in $X,Y$
is a priori in the quotient field $\bZ(a,b)$ of the polynomial ring $\bZ[a,b]$, 
and indeed 
\begin{gather*}
\Disc(\Phi_n^{**})\in
\begin{cases}
 \bZ[a,b]_{2d} & \text{if }n=1\\ 
 \bZ[a,b]_{2d_n(d_n-1)/(d-1)} & \text{if }n>1
\end{cases}
\end{gather*}
(see Remark \ref{th:integer} below).
For each $n\in\bN$ and each $\ell\in\bN$ (not necessarily dividing $n$ and) less than $n$,
we have the homogeneous resultant
\begin{gather*}
\Res(\Phi_n^{**},\Phi_{\ell}^{**})\in
 \begin{cases}
 \bZ[a,b]_{2d_nd/(d-1)} &\text{if }\ell=1\\
 \bZ[a,b]_{2d_nd_\ell/(d-1)} &\text{if }\ell>1
\end{cases}
\end{gather*}
between the homogeneous polynomials 
$\Phi_n^{**}$ and $\Phi_{\ell}^{**}$ in $X,Y$.

Recall the descriptions \eqref{eq:regular} of the ring $\bZ[\Rat_d]$
of regular functions on $\Rat_d$ and \eqref{eq:regularunit} 
of the unit group $(\bZ[\Rat_d])^*$ of it,
and note that $(d(d-1))\mid d_n$ if $n>1$ (by \eqref{eq:mobiusvanish} below).
For every $n\in\bN$ (resp.\ for every $n\in\bN$ and 
every $\ell\in\bN$ less than $n$), setting
\begin{gather*}
 N_d(n):=(2d)^{-1}\deg(\Disc(\Phi^{**}_{n}))\in\bN
\quad\bigl(\text{resp. }
N_d(n,\ell):=(2d)^{-1}\deg(\Res(\Phi_{n}^{**},\Phi_{\ell}^{**}))\in\bN
\bigr),
\end{gather*}
we have the regular function on $\Rat_d$ 
\begin{gather}
D_n^{**}:=\rho_d^{-N_d(n)}\Disc(\Phi^{**}_{n})\in\bZ[\Rat_d]
\quad\Bigl(\text{resp. }
R_{n,\ell}^{**}:=\rho_d^{-N_d(n,\ell)}\Res(\Phi_{n}^{**},\Phi_{\ell}^{**})\in\bZ[\Rat_d]\Bigr).
\label{eq:discresl}
\end{gather}
From \eqref{eq:formaldynam}, we have
$R_{n,\ell}^{**}\in(\bZ[\Rat_d])^*(=\bZ^*=\{\pm 1\})$ if $\ell\nmid n$ and $\ell<n$. 

For every $n\in\bN$, 
we call the hypersurface $V(D_n^{**})$ in $\Rat_d$ defined by (the zeros of)
$D_n^{**}\in\bZ[\Rat_d]$
the $n$-th formally exact parabolic bifurcation hypersurface in $\Rat_d$,
which contains all the hypersurfaces $V(R_{n,\ell}^{**})$ in 
$\Rat_d$ defined by (the zeros of) $R_{n,\ell}^{**}\in\bZ[\Rat_d]$
for every $\ell\in\bN$ dividing $n$ and less than $n$
(by \eqref{eq:formaldynam} and the chain rule),  
and call the union 
\begin{gather*}
 \bigcup_{n\in\bN}V(D_n^{**})
\end{gather*}
the parabolic bifurcation locus in $\Rat_d$.

Next, for every $n\in\bN$, there is a so called
$n$-th formally exact multiplier polynomial
\begin{gather*}
 p_n^{**}(T;[a:b])\in\bigl(\bZ[\Rat_d]\bigr)[T]
\end{gather*}
in $T$ of degree $d_n/n$
so that, letting $\Omega$ be an algebraically closed field 
of characteristic $0$ (e.g.\ $\bC$), for
every individual $f=f_{[a:b]}\in\Rat_d(\Omega)$, we have
\begin{gather}
 \bigl(p_n^{**}(T;f)\bigr)^n=\prod_{j=1}^{d_n}\bigl((f^n)'(z_{f,j}^{(n)})-T\bigr)\in\Omega[T]\label{eq:factormult}
\end{gather}
(the left hand side in which is 
not the $n$-th iteration but 
the $n$-th power of $p_n^{**}(T;f)$), where we write as $[\Phi_{F,n}^{**}=0]=\sum_{j=1}^{d_n}[z_{f,j}^{(n)}]$ 
the effective divisor 
on $\bP^1(\Omega)$ defined by the $d_n$ roots in $\bP^1(\Omega)$ 
of $\Phi_{F,n}^{**}\in\Omega[X,Y]_{d_n}$ for a lift $F$ (to $\Omega^2$) of $f$,
taking into account their multiplicities
(so the sequence $(z_{f,j}^{(n)})_{j=1}^{d_n}$ in $\bP^1(\Omega)$
is independent of the choice of $F$).
See \cite[Theorem 4.5 and the paragraph before that]{Silverman98} 
for a full account over $\bZ$, and also \cite[\S4.5]{SilvermanDynamics} 
over $\bQ$.

When $T=1$, we have $V(p_n^{**}(1;\,\cdot\,))=V(D_n^{**})$ (up to multiplicities),
where $V(p_n^{**}(1;\,\cdot\,))$ is
the hypersurface in $\Rat_d$ defined by 
(the zeros of) $p_n^{**}(1;\,\cdot\,)\in\bZ[\Rat_d]$.
On the other hand,
for every $n\in\bN$ and every $\ell\in\bN$ dividing $n$ and 
less than $n$, the (so called) $\frac{n}{\ell}$-th cyclic resultant for 
the polynomial $p_{\ell}^{**}(T;[a:b])$ in $T$ gives the regular function
\begin{gather*}
[a:b]\mapsto\Delta_{n,\ell}(f_{[a:b]}):=\Res\bigl(C_{\frac{n}{\ell}}(T),p_{\ell}^{**}(T;[a:b])\bigr)
\in\bZ[\Rat_d]
\end{gather*}
on $\Rat_d$; letting $\Omega$ be an algebraically
closed field of characteristic $0$, for every individual $f\in\Rat_d(\Omega)$, 
in terms of the above $(z_{f,j}^{(\ell)})_{j=1}^{d_\ell}$ in $\bP^1(\Omega)$,
we have 
\begin{gather}
 \bigl(\Delta_{n,\ell}(f)\bigr)^\ell=\prod_{j=1}^{d_\ell}C_{\frac{n}{\ell}}\bigl((f^\ell)'(z_{f,j}^{(\ell)})\bigr)
=\prod_{\omega:\,\text{an $\frac{n}{\ell}$-th primitive root of unity (in }\overline{\bQ})}\bigl(p_\ell^{**}(\omega;f)\bigr)^\ell.\label{eq:cyclicsubst}
\end{gather}
In particular, 
$V(\Delta_{n,\ell})=V(R_{n,\ell}^{**})$ (up to multiplicities), 
where $V(\Delta_{n,\ell})$ is
the hypersurface in $\Rat_d$ defined by (the zeros of) $\Delta_{n,\ell}\in\bZ[\Rat_d]$.
For every $n\in\bN$, we also set a rational function
\begin{gather*}
\Delta_{n,n}:=
\frac{p_n^{**}(1;\,\cdot\,)}{
\prod_{\ell:\,\mid n\text{ and }<n}\Delta_{n,\ell}}
\quad\text{on }\Rat_d.\label{eq:collide}
\end{gather*}
\begin{remark}
 Under our notational convention, it would be more consistent to write 
 $\Delta_{n,n},\Delta_{n,\ell}$ as $\Delta_{n,n}^{**},\Delta_{n,\ell}^{**}$,
 respectively, but for simplicity, we omit 
``**'' for them here and below.
\end{remark}

\subsection{Main results}
Our principal result is the following computation of 
both the defining regular function $D_n^{**}$ on $\Rat_d$
of the $n$-th formally exact parabolic bifurcation hypersurface 
$V(D_n^{**})$, $n\in\bN$,
and the regular function $R_{n,\ell}^{**}$ on $\Rat_d$,
$\ell \mid n$ and $<n$, in terms of $\Delta_{n,\ell}$
and $\Delta_{n,n}$.

\begin{mainth}\label{th:exact}
Fix an integer $d>1$. 
Then for every $n\in\bN$,
we have
\begin{gather*}
D_n^{**}=\pm(\Delta_{n,n})^n\cdot\prod_{\ell:\,\mid n\text{ and }<n}(\Delta_{n,\ell})^{n-\ell}\quad\text{on }\Rat_d,
\end{gather*}
and for every $n\in\bN$ and every $\ell\in\bN$ dividing $n$ and less than $n$, we have
\begin{gather*}
R_{n,\ell}^{**}=\pm(\Delta_{n,\ell})^\ell\quad\text{on }\Rat_d,
\end{gather*}
where by ``$\pm$'' we mean ``up to sign''.
Moreover, for every $n\in\bN$, we have $\Delta_{n,n}\in\bZ[\Rat_d]$.
\end{mainth}

The proof of Theorem \ref{th:exact} is based on (complex) 
potential theoretic computations
while the resulting formulas in Theorem \ref{th:exact} are algebraic.
For a (monic) polynomial, Theorem \ref{th:exact} recovers \cite[the first equality in Theorem A and Theorem 2.2]{MortonVivaldi95}.
As a consequence of Theorem \ref{th:exact},
we have the following geometric descriptions of both 
the $n$-th formally exact parabolic bifurcation hypersurface
$V(D_n^{**})$ in $\Rat_d$, $n\in\bN$,
and its subhypersurfaces $V(R_{n,\ell}^{**})$ in $\Rat_d$,
$\ell\mid n$ and $<n$, in terms of effective divisors in $\Rat_d$.

\begin{maincoro}\label{th:geometric}
Fixing $d>1$, for every $n\in\bN$, we have
\begin{gather*}
[D_n^{**}=0]
=n[\Delta_{n,n}=0]
+\sum_{\ell:\,\mid n\text{ and }<n}(n-\ell)[\Delta_{n,\ell}=0]
\quad\text{and}\\
[p_n^{**}(1;\,\cdot\,)=0]=[\Delta_{n,n}=0]+
\sum_{\ell:\,\mid n\text{ and }<n}[\Delta_{n,\ell}=0]
\quad\text{on }\Rat_d,
\end{gather*}
and for every $\ell\in\bN$ dividing $n$ and 
less than $n$,
$[R_{n,\ell}^{**}=0]=\ell[\Delta_{n,\ell}=0]$ on $\Rat_d$.
\end{maincoro}

From a viewpoint of 
dynamics, it would be more natural to study each $f$ up to a projective
coordinate changes of the source and target copies of $\bP^1$
rather than $f$ itself. Over $\bC$, the set $\rM_d(\bC)=\Rat_d(\bC)/\PSL(2,\bC)$
of all M\"obius transformation conjugacy classes $[f]$ 
of rational functions $f\in\Rat_d(\bC)$ is called the dynamical 
moduli space of degree $d$ rational functions on $\bP^1(\bC)$.
To equip this set $\rM_d(\bC)$ with a natural
complex affine algebraic variety structure,
we (abstractly) introduce the dynamical moduli space 
\begin{gather*}
 \rM_d:=\Rat_d//\SL_2 
\end{gather*}
of degree $d$ rational functions on $\bP^1$ as
the (GIT-)quotient of $\Rat_d$ by the (GIT-stable) conjugation action on $\Rat_d$
of the linear algebraic group scheme $\SL_2$(, and then
the variety (rather than just a set) $\rM_d(\bC)$ is obtained
by the base extension (from $\Spec\bZ$ to $\Spec\bC$) of the scheme $\rM_d$;
see \cite[\S 2]{Silverman98} over $\bZ$, and also 
\cite[\S 4.4]{SilvermanDynamics} over $\bQ$). 

The regular functions
$p_n^{**}(1;\,\cdot\,)$, $\Delta_{n,n}$, $\Delta_{n,\ell}$,
$D_n^{**}$, and
$R_{n,\ell}^{**}$
on $\Rat_d$ are all invariant under the action of $\SL_2$ on 
$\Rat_d$ (see Remark \ref{th:counting} below for the $\SL_2$-invariance of 
$D_n^{**}$, $n\in\bN$, and that
of $R_{n,\ell}^{**}$, $n>\ell$),  
and descend to $\rM_d$
as regular functions there.
Hence the assertions similar to those in Theorem \ref{th:exact}
and Corollary \ref{th:geometric} also hold 
when replacing $\Rat_d$ with $\rM_d$.

\section{Background}
\label{sec:background}
\subsection{Complex dynamics and potential theory}\label{sec:complexpotential}
Let 
\begin{gather*}
 \pi:\bC^2\setminus\{(0,0)\}\to\bP^1(=\bP^1(\bC)) 
\end{gather*}
be the canonical
projection so that $\pi(0,1)=\infty$ (and that $\pi(Z_0,Z_1)=Z_1/Z_0$ if $Z_0\neq 0$), and $\omega$ be the Fubini-Study
area element on $\bP^1$ normalized as $\omega(\bP^1)=1$. Set 
\begin{gather*}
 (Z_0,Z_1)\wedge(W_0,W_1)=Z_0W_1-Z_1W_0\in\bZ[Z_0,Z_1,W_0,W_1]_2,
\end{gather*}
so that for $Z,W\in\bC^2$ and $A\in\GL(2,\bC)$, 
\begin{gather}
 (AZ)\wedge(AW)=(\det A)(Z\wedge W),\label{eq:determinant} 
\end{gather}
and let $\|\cdot\|$ be the Euclidean norm
on $\bC^2$. The normalized chordal
metric on $\bP^1$ is defined as
\begin{gather*}
 [z,w]_{\bP^1}=\frac{|Z\wedge W|}{\|Z\|\cdot\|W\|}\quad\text{on }\bP^1\times\bP^1,
\end{gather*}
where $Z\in\pi^{-1}(z),W\in\pi^{-1}(w)$, so that $[0,\infty]_{\bP^1}=1$.
The (generalized) 
Laplacian $\Delta$ on $\bP^1$ is normalized so that
$\Delta(\log[\cdot,\infty]_{\bP^1})=\delta_\infty-\omega$ on $\bP^1$,
where $\delta_z$ is the Dirac measure on $\bP^1$ at $z\in\bP^1$.

Pick $f\in\Rat_d(\bC)$, $d>1$, and a lift $F\in(\bC[X,Y]_d)^2$ 
(to $\bC^2$) of $f$,
which is unique up to multiplication in $\bC^*=\bC\setminus\{0\}$.
Then for every $n\in\bN$ and every $c\in\bC^*$,
$(cF)^n=c^{(d^n-1)/(d-1)}F^n\in(\bC[X,Y]_{d^n})^2$ is a lift (to $\bC^2$) of $f^n$.
From the homogeneity of $F$, the uniform limit
\begin{gather*}
 G^F:=\lim_{n\to\infty}\frac{\log\|F^n\|}{d^n}\quad\text{on }\bC^2\setminus\{(0,0)\}
\quad(F^n=F\circ\cdots\circ F,\ n\text{ times})
\end{gather*}
exists and is called the escaping rate function of $F$ on $\bC^2\setminus\{(0,0)\}$. 
We note that
for every $c\in\bC^*$, $G^F(cZ)=G^F(Z)+\log|c|$ 
and $G^{cF}=G^F+(\log|c|)/(d-1)$ on $\bC^2\setminus\{(0,0)\}$, that
\begin{gather}
 G^F\circ F=d\cdot G^F\quad\text{on }\bC^2\setminus\{(0,0)\},\label{eq:pullbackGreen}
\end{gather} 
and that for every $n\in\bN$, 
$G^{F^n}=G^F$ on $\bC^2\setminus\{(0,0)\}$.

The continuous function
\begin{gather*}
 g_F:=G^F-\log\|\cdot\|\quad\text{on }\bP^1
\end{gather*}
is called the dynamical Green function of $F$ on $\bP^1$,
and the probability Radon measure $\mu_f:=\Delta g_F+\omega$ on $\bP^1$
is independent of the choice of $F$ and is called the equilibrium (or canonical) measure of $f$
on $\bP^1$. The homogeneous resultant of a lift (to $\bC^2$)
$F=F(X,Y;a,b)=(F_0,F_1)$ of $f_{[a:b]}\in\Rat_d(\bC)$ 
($a=(a_0,\ldots,a_d),b=(b_0,\ldots,b_d)$ for short) is
\begin{gather*}
 \Res F:=\Res(F_0,F_1)=\rho_d(a,b)\in\bC^*,
\end{gather*}
where $\Res(F_0,F_1)$ is the homogeneous resultant
between $F_0,F_1\in\bC[X,Y]_d$ (see Subsection \ref{sec:resultant} below).
For every $c\in\bC^*$, we also have $\Res(cF)=c^{2d}\Res F$,
so that the function
\begin{gather*}
 g_f:=g_F-\frac{\log|\Res F|}{2d(d-1)}\quad\text{on }\bP^1
\end{gather*}
is independent of the choice of $F$ and is called the dynamical Green function
of $f$ on $\bP^1$. 
We also note that for every $n\in\bN$, $g_{f^n}=g_f$ on $\bP^1$ or equivalently
$|\Res(F^n)|=|\Res F|^{(d^n(d^n-1))/(d(d-1))}$. 

For each lift $F$ (to $\bC^2$) of $f$,
the $g_F$-kernel function on $\bP^1$ is defined by the function
$\Phi_{g_F}(z,w):=\log[z,w]_{\bP^1}-g_F(z)-g_F(w)$
on $\bP^1\times\bP^1$ (there would be no confusions between this and dynatomic polynomials, 
in notations), so that $\Delta(\Phi_{g_F}(\cdot,w))=\delta_w-\mu_f$
on $\bP^1$ for every $w\in\bP^1$. We have 
$\int_{\bP^1}\Phi_{g_f}(\cdot,w)\mu_f(w)\equiv-(\log|\Res F|)/(d(d-1))$
on $\bP^1$
since $\Delta\int_{\bP^1}\Phi_{g_f}(\cdot,w)\mu_f(w)=\mu_f-\mu_f=0$ on $\bP^1$
(by the Fubini theorem)
and moreover the energy/capacity formula
\begin{gather*}
 \int_{\bP^1\times\bP^1}\Phi_{g_F}(\mu_f\times\mu_f)=-\frac{\log|\Res F|}{d(d-1)} 
\end{gather*}
also holds (due to DeMarco \cite{DeMarco03};
for a simple proof which also works in non-archimedean setting,
see \cite[Appendix A]{Baker09} or \cite[Appendix]{OS11})). 
In particular, 
the $g_f$-kernel function
\begin{gather*}
 \Phi_{g_f}(z,w):=\log[z,w]_{\bP^1}-g_f(z)-g_f(w)\quad\text{on }\bP^1\times\bP^1
\end{gather*}
satisfies
$\int_{\bP^1}\Phi_{g_f}(\cdot,w)\mu_f(w)\equiv 0$ on $\bP^1$.

The following Riesz decomposition-type formula
\begin{gather}
 \Phi_{g_f}(f,\Id_{\bP^1})=\int_{\bP^1}\Phi_{g_f}(\cdot,w)[f=\Id_{\bP^1}](w)\quad\text{on }\bP^1\label{th:byparts}
\end{gather}
(\cite[Lemma 4.2]{OS15}, the most involved point is that
the harmonic part $\int_{\bP^1}\Phi_{g_f}(f,\Id_{\bP^1})\mu_f\in\bR$ 
of $\Phi_{g_f}(f,\Id_{\bP^1})$ also vanishes on $\bP^1$) 
plays a key role. Here, writing
\begin{gather}
 F(Z)\wedge Z=\prod_{j=1}^{d+1}(Z\wedge q_{F,j})\quad\text{on }\bC^2\label{eq:factorfix}
\end{gather}
for some $q_{F,1},\ldots,q_{F,d+1}\in\bC^2\setminus\{(0,0)\}$,
the sequence $(\pi(q_{F,j}))_{j=1}^{d+1}$ is independent
of the choice of $F$, up to permutation, and
the effective divisor $[f=\Id_{\bP^1}]=\sum_{j=1}^{d+1}[\pi(q_{F,j})]$ on $\bP^1$ 
defined by the solutions in $\bP^1$ of the equation $f=\Id_{\bP^1}$
taking into account their multiplicities is regarded as 
a positive Radon measure 
\begin{gather*}
[f=\Id_{\bP^1}]=\sum_{j=1}^{d+1}\delta_{\pi(q_{F,j})}\quad\text{on }\bP^1.
\end{gather*}
For every $z\in\Fix(f)$,
(denoting the multiplier of the cycle $z$ of $f$ (having the period $1$) by $f'(z)$
even when $z=\infty$, for simplicity,) 
if this $z$ is simple in that $f'(z)\neq 1$, then we have
\begin{gather}
 \lim_{\bP^1\ni w\to z}\bigl(\Phi_{g_f}(f(w),w)-\Phi_{g_f}(z,w)\bigr)
=\lim_{\bP^1\ni w\to z}\log\frac{[f(w),w]_{\bP^1}}{[z,w]_{\bP^1}}=\log|f'(z)-1|,\label{eq:diffderiv}
\end{gather}
which together with \eqref{th:byparts} yields
\begin{gather}
 \log|f'(z)-1|=\int_{\bP^1\setminus\{z\}}\Phi_{g_f}(z,w)[f=\Id_{\bP^1}](w).\label{eq:multiplier}
\end{gather}

In the rest of this subsection, let us assume that $F$ satisfies $|\Res F|=1$,
so that
\begin{gather}
\Phi_{g_f}\bigl(\pi(Z),\pi(W)\bigr)=\log|Z\wedge W|-G^F(Z)-G^F(W) 
\quad\text{on }\bigl(\bC^2\setminus\{(0,0)\}\bigr)^2.\label{eq:homogkernel} 
\end{gather}
Another consequence of \eqref{th:byparts}(, \eqref{eq:homogkernel},
the factorization \eqref{eq:factorfix} of $F\wedge\Id_{\bC^2}$,
and \eqref{eq:pullbackGreen}) is the vanishing
\begin{gather}
 \sum_{j=1}^{d+1}G^F(q_{F,j})=0.\label{eq:Greensum}
\end{gather}
Hence for every $n\in\bN$, noting that
$F^n\wedge\Id_{\bC^2}=\prod_{m\,\mid n}\Phi_{F,m}^{**}$
(by M\"obius inversion)
and writing each $\Phi_{F,m}^{**}$ as
$\Phi_{F,m}^{**}(Z)=\prod_{j=1}^{d_m}(Z\wedge Z_{F,j}^{(m)})$
for some $Z_{F,1}^{(m)},\ldots,Z_{F,d_m}^{(m)}\in\bC^2\setminus\{(0,0)\}$,
by the vanishing \eqref{eq:Greensum} applied to $F^n$ 
and by M\"obius inversion,
we also have
\begin{gather}
 \sum_{j=1}^{d_n}G^F(Z_{F,j}^{(n)})=0,\label{eq:Greenexact}
\end{gather}
which will be frequently used in this paper.

\subsection{M\"obius function and inversions}\label{th:mobius}
For the details, see e.g.\ \cite[\S2.7]{Apostol}.
The M\"obius function $\mu:\bN\to\{0,\pm 1\}$ is defined so that
$\mu(n)=1\cdot(-1)^N$ if 
$n$ is a product of some $N$
distinct positive prime numbers, $N\in\bN\cup\{0\}$, (so $\mu(1)=1$) 
and $\mu(n)=0$ in other cases. 
We repeatedly use
\begin{gather}
  \sum_{m\mid n}\mu(m)=0\quad\text{if }n>1\label{eq:mobiusvanish} 
\end{gather}
in this paper. For any sequences 
 $(A_k)_{k\in\bN},(B_\ell)_{\ell\in\bN}$ in a commutative group,
 \begin{gather*}
 A_k=\prod_{\ell\mid k}B_\ell\text{ for every }k\in\bN\quad\text{if and only if}\quad
 B_\ell=\prod_{k\mid\ell}(A_k)^{\mu(\frac{\ell}{k})}\text{ for every }\ell\in\bN,
 \end{gather*}
 and then we say the sequences $(A_k),(B_\ell)$ are obtained by (multiplicative)
 M\"obius inversions of each other. The same thing could also be
stated additively (in an abelian group).

\subsection{Homogeneous resultant and discriminant}\label{sec:resultant}
For any homogeneous polynomials $P=\sum_{i=0}^da_iX^{d-i}Y^i\in R[X,Y]_d$ and $Q=\sum_{j=0}^eb_jX^{e-j}Y^j\in R[X,Y]_e$ over an integral domain $R$, 
the homogeneous resultant between $P$ and $Q$ is defined as
\begin{gather*}
\Res(P,Q):=
\det\begin{pmatrix}
     a_0 & \cdots & a_{e-1} & \cdots & a_{d-1} &a_d  &  & \\
     \text{\huge{0}} & \ddots & \vdots  &  &  &  & \ddots & \text{\huge{0}}   \\
     & & a_0     & a_1  & \cdots & \cdots & a_{d-1} & a_d\\
     b_0 & \cdots & b_{e-1} & b_e  \\
     \text{\huge{0}}  & \ddots & &  & \vdots  &  & \ddots & \text{\huge{0}}  \\
       &  & &  & b_0 &  \cdots & b_{e-1} & b_e 
    \end{pmatrix}\in R,
\end{gather*}
so that $\Res(P,Q)=(-1)^{de}\Res(Q,P)$ and that
$\Res(cP,Q)=c^e\Res(P,Q)$ for every $c\in R$.

Letting $F$ be the quotient field of $R$,
any $P\in R[X,Y]_d$ factors as 
\begin{gather*}
 P(X,Y)=\prod_{i=1}^d\bigl((X,Y)\wedge A_i\bigr) 
\end{gather*}
over an algebraic closure
$\overline{F}$ of $F$, where
$A_i\in\overline{F}^2\setminus\{(0,0)\}$ for each $i\in\{1,\ldots,d\}$.
Then for every $Q\in R[X,Y]_d$,
\begin{gather}
 \Res(P,Q)=\prod_{i=1}^d Q(A_i),\label{eq:resultantsubst}
\end{gather}
and the homogeneous discriminant of $P$ is defined as
\begin{gather}
 \Disc(P):=\prod_{i=1}^d\,\prod_{\ell\in\{1,\ldots,d\}\setminus\{i\}}(A_i\wedge A_\ell)\in F,\label{eq:homogdiscr}
\end{gather}
so that $\Disc(cP)=c^{2(d-1)}\Disc(P)$ for every $c\in R$.
For more details, see e.g.\ \cite[\S2.4]{SilvermanDynamics}.

\begin{remark}\label{th:counting}
In terms of the notations 
in Subsections \ref{sec:complexpotential} and \ref{sec:dynatomic},
for every $A=(a_{st})\in\GL(2,\bC)$, 
regarding $A$ as a linear automorphism $(X,Y)\mapsto(a_{11}X+a_{12}Y,a_{21}X+a_{22}Y)$
of $\bC^2$ and setting a M\"obius transformation
$a(z):=(a_{21}+a_{22}z)/(a_{11}+a_{12}z)$ of $\bP^1(\bC)$,
the conjugation $A\circ F\circ A^{-1}\in(\bC[X,Y]_d)^2$ of $F$
is a lift (to $\bC^2$) of the conjugation $a\circ f\circ a^{-1}\in\Rat_d(\bC)$ of $f$,
and using \eqref{eq:determinant} repeatedly, we have 
\begin{gather*}
 ((A\circ F\circ A^{-1})(Z))\wedge Z=(\det A)^{-d}\prod_{j=1}^{d+1}(Z\wedge(Aq_{F,j})) 
\end{gather*}
on $\bC^2$.
This together with 
\eqref{eq:determinant}, \eqref{eq:resultantsubst}, \eqref{eq:homogdiscr}, and
\begin{gather*}
 \Res(A\circ F\circ A^{-1})=(\det A)^{d-d^2}\Res F
\end{gather*}
for every $A\in\GL(2,\bC)$ (see, e.g., \cite[Exercise 2.12]{SilvermanDynamics})
concludes the $\SL_2$-invariance of 
$D_n^{**}$, $n\in\bN$, and that of 
$R_{n,\ell}^{**}$,
$n>\ell$, on $\Rat_d$.
\end{remark}

\subsection{A proposition from algebra}
For completeness, we include a proof of the following.
\begin{proposition}\label{th:unit}
 Let $R$ be a uniquely factorization domain $($UFD$)$, $N\in\bN$, and $F$ be an irreducible
 homogeneous polynomial in $R[X_0,\ldots,X_N]$ of degree $>0$, and set
 $U:=\bP^N_R\setminus V(F)$, where $V(F)$ be the hypersurface in $\bP^N_R$
 defined by $($the zeros in $\bP^N_R$ of$)$ $F$. 
 Then the unit group $(R[U])^*$ of the ring $R[U]$ of regular functions
 on $U$ equals $R^*$.
\end{proposition}

\begin{proof}
 Pick $\phi\in(R[U])^*$. Then $\phi=P/F^{(\deg P)/(\deg F)}$ 
 and $\phi^{-1}=Q/F^{(\deg Q)/(\deg F)}$ for some
 homogeneous $P,Q\in R[X_0,\ldots,X_N]$ such that
 $(\deg P)/(\deg F),(\deg Q)/(\deg F)\in\bN\cup\{0\}$.
 Since $R[X_0,\ldots,X_N]$ is also UFD,
 $PQ=F^{(\deg P)/(\deg F)+(\deg Q)/(\deg F)}$ implies
 $P=rF^m$ for some $r\in(R[X_0,\ldots,X_n])^*=R^*$ and some $m\in\bN\cup\{0\}$,
 so $\phi=(rF^m)/F^{(m\deg F)/\deg F}=r\in R^*$.
\end{proof}

\section{Proof of Theorem \ref{th:exact}}
\label{sec:proofself}
It is enough to show Theorem \ref{th:exact} working over $\bC$. 
Fix $d>1$, and pick $n\in\bN$. 
For every individual (generic) $f\in\Rat_d(\bC)$
and every lift $F$ (to $\bC^2$) of $f$, 
as in Section \ref{sec:complexpotential},
we have a factorization
$F^n(Z)\wedge Z=\prod_{j=1}^{d^n+1}(Z\wedge q_{F,j}^{(n)})$ 
on $\bC^2$ of $F^n\wedge\Id_{\bC^2}$, 
where $q_{F,j}^{(n)}\in\bC^2\setminus\{(0,0)\}$
for every $j\in\{1,\ldots,d^n+1\}$.
Moreover, for every $\ell\in\bN$, we also have a factorization
\begin{gather*}
\Phi_{F,\ell}^{**}(Z)=\prod_{j=1}^{d_\ell}(Z\wedge Z_{F,j}^{(\ell)}) \end{gather*}
on $\bC^2$ of $\Phi_{F,\ell}^{**}$, where $Z_{F,j}^{(\ell)}\in\bC^2\setminus\{(0,0)\}$ for every $j\in\{1,\ldots,d_\ell\}$, 
and setting 
\begin{gather}
 z_{f,j}^{(\ell)}:=\pi\bigl(Z_{F,j}^{(\ell)}\bigr)\in\bP^1=\bP^1(\bC)\quad
\text{for each }j\in\{1,\ldots,d_\ell\},\label{eq:multiset} 
\end{gather}
the sequence $(z_{f,j}^{(\ell)})_{j=1}^{d_\ell}$ is 
independent of the choice of $F$, up to permutation, so that
$[\Phi_{F,\ell}^{**}=0]=\sum_{j=1}^{d_\ell}[z_{f,j}^{(\ell)}]$ 
(as an effective divisor) on $\bP^1$.

Without loss of generality, 
we normalize the labeling of the sequence $(z_{f,j}^{(n)})_{j=1}^{d_n}$
so that for every $r\in\{0,1,\ldots,\frac{d_n}{n}-1\}$ 
and every $k\in\{1,\ldots,n-1\}$, 
\begin{gather}
f^k(z_{f,1+nr}^{(n)})=z_{f,1+nr+k}^{(n)}
\quad\Bigl(\text{and}\quad f^n(z_{f,1+nr}^{(n)})=z_{f,1+nr}^{(n)}\Bigr),\label{eq:labelingformal}
\end{gather}
and then for every $r\in\{0,1,\ldots,\frac{d_n}{n}-1\}$, set
the homogeneous polynomial
\begin{gather}
 \Lambda_{F,r}^{(n)}(X,Y):=\prod_{k=0}^{n-1}((X,Y)\wedge Z_{F,1+nr+k}^{(n)})\in\bC[X,Y]_n
\quad\Bigl(\text{so that }\Phi_{F,n}^{**}=\prod_{r=0}^{\frac{d_n}{n}-1}\Lambda_{F,r}^{(n)}\Bigr).\label{eq:cyclemutual}
\end{gather}

The equalities \eqref{eq:collideprecise} and \eqref{eq:factoredmult}
below are substantial.

\subsection{}
We first claim that
for every $\ell\in\bN$ dividing $n$ and less than $n$,
\begin{gather}
\biggl|\frac{(\Res F)^{-N_d(n,\ell)}\Res(\Phi_{F,n}^{**},\Phi_{F,\ell}^{**})}{(\Delta_{n,\ell}(f))^\ell}\biggr|=1, 
\label{eq:collideprecise}
\end{gather} 
which yields $R_{n,\ell}^{**}/(\Delta_{n,\ell})^\ell\in(\bZ[\Rat_d])^*(=\bZ^*=\{\pm 1\})$
recalling \eqref{eq:regularunit} and
the definition of $R_{n,\ell}^{**}$ in \eqref{eq:discresl};
for, we note that
\begin{gather*}
\Phi_{F,n}^{**}(Z)
\Bigl(=\prod_{m\,\mid n}(F^m(Z)\wedge Z)^{\mu(\frac{n}{m})}\Bigr)
=\prod_{m:\,\ell\nmid m\mid n}(F^m(Z)\wedge Z)^{\mu(\frac{n}{m})}\times
\prod_{m:\,\ell\mid m\mid n}(F^m(Z)\wedge Z)^{\mu(\frac{n}{m})}
\end{gather*}
on $\bC^2$, and compute
\begin{align*}
&\prod_{m:\,\ell\nmid m\mid n}(F^m(Z)\wedge Z)^{\mu(\frac{n}{m})}
=\prod_{m:\,\ell\nmid m\mid n}\Bigl(\prod_{m'\,\mid m}\Phi_{F,m'}^{**}(Z)\Bigr)^{\mu(\frac{n}{m})}\\
=&
\prod_{m':\,\ell\nmid m'\mid n}\Bigl(\prod_{m:\,m'\mid m\mid n\text{ and }\ell\nmid m}\Phi_{F,m'}^{**}(Z)^{\mu(\frac{n}{m})}\Bigr)\\
=&\prod_{m':\,m'\nmid\ell\nmid m'\mid n}\Bigl(\prod_{m:\,m'\mid m\mid n\text{ and }\ell\nmid m}\Phi_{F,m'}^{**}(Z)^{\mu(\frac{n}{m})}\Bigr)\times
\prod_{m':\,\mid\ell\text{ and }<\ell}\Bigl(\prod_{m:\,m'\mid m\mid n\text{ and }\ell\nmid m}\Phi_{F,m'}^{**}(Z)^{\mu(\frac{n}{m})}\Bigr),
\end{align*}
where the first equality is from
$F^m(Z)\wedge Z=\prod_{m'\,\mid m}\Phi_{F,m'}^{**}(Z)$
by M\"obius inversion. Moreover, we have 
\begin{multline*}
\prod_{m':\,\mid\ell\text{ and }<\ell}\Bigl(\prod_{m:\,m'\mid m\mid n\text{ and }\ell\nmid m}\Phi_{F,m'}^{**}(Z)^{\mu(\frac{n}{m})}\Bigr)
=\prod_{m':\,\mid\ell\text{ and }<\ell}\Bigl(\prod_{m'':\,\frac{\ell}{m'}\nmid m''\mid\frac{n}{m'}}\Phi_{F,m'}^{**}(Z)^{\mu(\frac{n/m'}{m''})}\Bigr)\\
=\prod_{m':\,\mid\ell\text{ and }<\ell}\Phi_{F,m'}^{**}(Z)^{\sum_{m''\,\mid\frac{n}{m'}}\mu(\frac{n/m'}{m''})-\sum_{m'':\,\frac{\ell}{m'}\mid m''\mid\frac{n}{m'}}\mu(\frac{n/m'}{m''})}\equiv 1,
\end{multline*}
noting that by \eqref{eq:mobiusvanish},
\begin{gather*}
 \begin{cases}
\sum_{m''\,\mid\frac{n}{m'}}\mu\bigl(\frac{n/m'}{m''}\bigr)=0\quad\text{(also by }\frac{n}{m'}>1)\quad\text{and}\\
\sum_{m'':\,\frac{\ell}{m'}\mid m''\mid\frac{n}{m'}}\mu\bigl(\frac{n/m'}{m''}\bigr)
 =\sum_{m'''\,\mid\frac{n}{\ell}}\mu\bigl(\frac{n/\ell}{m'''}\bigr)=0
 \quad\text{(also by }\frac{n}{\ell}>1) 
\end{cases}
\end{gather*}
when $m'\in\bN$ divides $\ell$ and is $<\ell$
(this purely algebraic computation is similar to that in \cite[Proof of Theorem 2.2]{MortonVivaldi95}).

On the other hand, for every $m\in\bN$, assuming $|\Res F|=1$, we compute 
\begin{align*}
&\log\biggl|\prod_{m:\,\ell\mid m\mid n}(F^m(Z)\wedge Z)^{\mu(\frac{n}{m})}\biggr|
=\sum_{k\,\mid\frac{n}{\ell}}\mu\Bigl(\frac{n/\ell}{k}\Bigr)\log|F^{k\ell}(Z)\wedge Z|\\
=&\sum_{k\,\mid\frac{n}{\ell}}\mu\Bigl(\frac{n/\ell}{k}\Bigr)
\Bigl(\Phi_{g_f}\bigl(f^{k\ell}(\pi(Z)),\pi(Z)\bigr)+G^F(F^{k\ell}(Z))+G^F(Z)\Bigr)
\quad(\text{by }\eqref{eq:homogkernel})\\
=&\sum_{k\,\mid\frac{n}{\ell}}\mu\Bigl(\frac{n/\ell}{k}\Bigr)
\Bigl(\Phi_{g_f}\bigl(f^{k\ell}(\pi(Z)),\pi(Z)\bigr)+\bigl((d^{\ell})^k+1\bigr)\cdot G^F(Z)\Bigr)
\quad(\text{by }\eqref{eq:pullbackGreen})\\
=&\sum_{k\,\mid\frac{n}{\ell}}\mu\Bigl(\frac{n/\ell}{k}\Bigr)
\Phi_{g_f}\bigl(f^{k\ell}(\pi(Z)),\pi(Z)\bigr)
+(d^{\ell})_{\frac{n}{\ell}}\cdot G^F(Z)
\end{align*}
generically on $\bC^2$, where $(d^{\ell})_{\frac{n}{\ell}}
=\sum_{k\,\mid\frac{n}{\ell}}\mu(\frac{n/\ell}{k})((d^{\ell})^k+1)$;
recalling the definition \eqref{eq:multiset} of the sequence
$(z_{f,j}^{(m)})_{j=1}^{d_m}$, for every $j\in\{1,\ldots,d_\ell\}$,
we also have
\begin{align*}
&\lim_{\bC^2\ni Z\to Z_j^{(\ell)}}\sum_{k\,\mid\frac{n}{\ell}}\mu\Bigl(\frac{n/\ell}{k}\Bigr)\Phi_{g_f}\bigl(f^{k\ell}(\pi(Z)),\pi(Z)\bigr)\\
=&\lim_{\bC^2\ni Z\to Z_j^{(\ell)}}
\sum_{k\,\mid\frac{n}{\ell}}\mu\Bigl(\frac{n/\ell}{k}\Bigr)
\Bigl(\Phi_{g_f}\bigl(f^{k\ell}(\pi(Z)),\pi(Z)\bigr)-\Phi_{g_f}\bigl(\pi(Z),\pi(Z_j^{(\ell)})\bigr)\Bigr)\\
=&\sum_{k\,\mid\frac{n}{\ell}}\mu\Bigl(\frac{n/\ell}{k}\Bigr)
\log\bigl|(f^{k\ell})'(z_{f,j}^{(\ell)})-1\bigr|
\quad(\text{by }\eqref{eq:diffderiv})\\
=&\sum_{k\,\mid\frac{n}{\ell}}\mu\Bigl(\frac{n/\ell}{k}\Bigr)
\log\bigl|\bigl((f^{\ell})'(z_{f,j}^{(\ell)})\bigr)^k-1\bigr|
=\log\bigl|C_{\frac{n}{\ell}}\bigl((f^{\ell})'(z_{f,j}^{(\ell)})\bigr)\bigr|
\quad(\text{using the chain rule}),
\end{align*}
where the second equality follows from
$\sum_{m\,\mid\frac{n}{\ell}}\mu(\frac{n/\ell}{m})=0$ 
also by 
\eqref{eq:mobiusvanish} and $n/\ell>1$.

From the above computations, when $|\Res F|=1$, using \eqref{eq:resultantsubst} repeatedly,
we have
\begin{multline*}
\log\bigl|\Res(\Phi_{F,n}^{**},\Phi_{F,\ell}^{**})\bigr|\\
\biggl(=\sum_{j=1}^{d_\ell}\log\bigl|\Phi_{F,n}^{**}(Z_{F,j}^{(\ell)})\bigr|
=\log\biggl|\prod_{m':\,m'\nmid\ell\nmid m'\mid n}\Bigl(\prod_{m:\,m'\mid m\mid n\text{ and }\ell\nmid m}\Res(\Phi_{F,m'}^{**},\Phi_{F,\ell}^{**})^{\mu(\frac{n}{m})}\Bigr)\biggr|+0+\\
+\log\biggl|\prod_{j=1}^{d_\ell}C_{\frac{n}{\ell}}\bigl((f^{\ell})'(z_{f,j}^{(\ell)})\bigr)\biggr|
+(d^{\ell})_{\frac{n}{\ell}}\sum_{j=1}^{d_\ell}G^F(Z_{F,j}^{(\ell)})=\biggr)\\
=\log\prod_{m':\,m'\nmid\ell\nmid m'\mid n}\Bigl(\prod_{m:\,m'\mid m\mid n\text{ and }\ell\nmid m}|\Res(\Phi_{F,m'}^{**},\Phi_{F,\ell}^{**})|^{\mu(\frac{n}{m})}\Bigr)+\log\bigl|\bigl(\Delta_{n,\ell}(f)\bigr)^\ell\bigr|+0\quad\text{(by \eqref{eq:cyclicsubst}, \eqref{eq:Greenexact})},
\end{multline*}
so that, no matter whether $|\Res F|=1$, we have
\begin{multline*}
\biggl|\frac{(\Res F)^{-N_d(n,\ell)}\Res(\Phi_{F,n}^{**},\Phi_{F,\ell}^{**})}{(\Delta_{n,\ell}(f))^\ell}\biggr|\\
=\prod_{m':\,m'\nmid\ell\nmid m'\mid n}\Bigl(\prod_{m:\,m'\mid m\mid n\text{ and }\ell\nmid m}\Bigl|(\Res F)^{-N_d(m',\ell)}\Res(\Phi_{F,m'}^{**},\Phi_{F,\ell}^{**})\Bigr|^{\mu(\frac{n}{m})}\Bigr)=1
\end{multline*}
since if $m'\nmid\ell\nmid m'$, then
$(\Res F)^{-N_d(m',\ell)}\Res(\Phi_{F,m'}^{**},\Phi_{F,\ell}^{**})\in(\bZ[\Rat_d])^*(=\bZ^*=\{\pm 1\})$.
Hence the desired \eqref{eq:collideprecise} holds.

\subsection{}
Next, we claim that 
\begin{gather}
\bigl|\Delta_{n,n}(f)\bigr|^n
=\biggl|\frac{(\Res F)^{-N_d(n)}\Disc(\Phi_{F,n}^{**})}{\prod_{\ell:\,\mid n\text{ and }<n}(\Delta_{n,\ell}(f))^{n-\ell}}\biggr|=
\prod_{r=0}^{\frac{d_n}{n}-1}\prod_{s\neq r}
\Bigl|(\Res F)^{-M_d(n)}\Res(\Lambda_{F,r}^{(n)},\Lambda_{F,s}^{(n)})\Bigr|,
\label{eq:factoredmult}
\end{gather}
the first equality in which yields 
$D_n^{**}/((\Delta_{n,n})^n\cdot\prod_{\ell:\,\mid n\text{ and }<n}(\Delta_{n,\ell})^{n-\ell})\in(\bZ[\Rat_d])^*(=\bZ^*=\{\pm 1\})$
recalling \eqref{eq:regularunit} and
the definition of $D_n^{**}$ in \eqref{eq:discresl}, and in the second one in which,
setting
\begin{gather*}
 M_d(n):=\begin{cases}
	       1 & \text{if }n=1\\
	       d_n(d_n-n)/(d(d-1)) & \text{if }n>1
	      \end{cases}\in\bN,
\end{gather*}
the product $(\Res F)^{-M_d(n)}\prod_{r=0}^{\frac{d_n}{n}-1}\prod_{s\neq r}
\Res(\Lambda_{F,r}^{(n)},\Lambda_{F,s}^{(n)})$ is independent of the choice of $F$
from $\Res(cF)=c^{2d}\Res F$ and $(cF)^n=c^{(d^n-1)/(d-1)}\Res F$
for every $c\in\bC^*$;

{\bfseries (i)} for, by $F^n(Z)\wedge Z=\prod_{\ell\mid n}\Phi_{F,\ell}^{**}(Z)$ 
(by M\"obius inversion) and
the definition \eqref{eq:homogdiscr} of the homogeneous discriminants, we have
\begin{align*}
 \Disc(F^n\wedge\Id_{\bC^2})
=&\prod_{\ell|n}\prod_{j=1}^{d_\ell}\biggl(\prod_{\ell':|n\text{ and }\neq\ell}\,\prod_{k=1}^{d_{\ell'}}(Z_{F,j}^{(\ell)}\wedge Z_{F,k}^{(\ell')})\times\prod_{k\in\{1,\ldots,d_\ell\}\setminus\{j\}}(Z_{F,j}^{(\ell)}\wedge Z_{F,k}^{(\ell)})\biggr)\\
=&\prod_{\ell|n}\,\prod_{\ell':|n\text{ and }\neq\ell}\,\prod_{j=1}^{d_\ell}\prod_{k=1}^{d_{\ell'}}(Z_{F,j}^{(\ell)}\wedge Z_{F,k}^{(\ell')})\times\prod_{\ell|n}\Disc(\Phi_{F,\ell}^{**})\\
=&
\prod_{(\ell,\ell')\in\bN^2:\,\ell\mid n,\,\ell'\mid n,\text{ and }\ell\neq\ell'}\Res(\Phi_{F,\ell}^{**},\Phi_{F,\ell'}^{**})\times \prod_{\ell\,\mid n}\Disc(\Phi_{F,\ell}^{**}),
\end{align*}
using \eqref{eq:resultantsubst} repeatedly.
Assume now that $|\Res F|=1$. Then we also compute not only
\begin{align*}
&\log\bigl|\Disc(F^n\wedge\Id_{\bC^2})\bigr|\\
=&\sum_{j=1}^{d^n+1}\sum_{k\in\{1,\ldots,d^n+1\}\setminus\{j\}}\Phi_{g_f}\bigl(\pi(q_{F,j}^{(n)}),\pi(q_{F,k}^{(n)})\bigr)+(d^n+1)\sum_{j=1}^{d^n+1}G^F(q_{F,j}^{(n)})
\quad(\text{by }\eqref{eq:homogkernel})\\
=&\sum_{j=1}^{d^n+1}\log|(f^n)'(\pi(q_{F,j}^{(n)}))-1|+0
\quad(\text{by }\eqref{eq:multiplier}\text{ and }\eqref{eq:Greensum})\\
=&\sum_{\ell\,\mid n}\sum_{j=1}^{d_\ell}\log\bigl|\bigl((f^\ell)'(z_{f,j}^{(\ell)})\bigr)^{\frac{n}{\ell}}-1\bigr|
\quad\bigl(\text{by }F^n(Z)\wedge Z=\prod_{\ell\,\mid n}\Phi_{F,\ell}^{**}(Z)\text{ and the chain rule}\bigr)\\
=&\sum_{\ell\,\mid n}\sum_{j=1}^{d_\ell}
\log\biggl|\prod_{\ell'\,\mid \frac{n}{\ell}}
C_{\ell'}\bigl((f^\ell)'(z_{f,j}^{(\ell)})\bigr)\biggr|
\quad(\text{by M\"obius inversion of the sequence }(C_i)_i)\\
=&\sum_{\ell\,\mid n}\,\sum_{\ell'\,\mid \frac{n}{\ell}}
\log\biggl|\prod_{j=1}^{d_\ell}
C_{\ell'}\bigl((f^\ell)'(z_{f,j}^{(\ell)})\bigr)\biggr|
=\sum_{\ell\,\mid n}\Bigl(\log\bigl|p_{\ell}^{**}(1;f)\bigr|^\ell
+\sum_{\ell':\,1<\ell'\mid\frac{n}{\ell}}\log\bigl|\Delta_{(\ell'\ell),\ell}(f)\bigr|^\ell\Bigr)\\
=&\sum_{\ell\,\mid n}\log\bigl|p_{\ell}^{**}(1;f)\bigr|^\ell
+\sum_{(\ell,\ell'')\in\bN^2:\,\ell|\ell''|n\text{ and }\ell<\ell''}\log\bigl|\Delta_{\ell'',\ell}(f)\bigr|^\ell
\end{align*}
(where the 6th equality is by 
the definition
\eqref{eq:factormult} of $p_{\ell}^{**}$ and
the first equality in \eqref{eq:cyclicsubst}) but also
\begin{multline*}
 \prod_{(\ell,\ell')\in\bN^2:\,\ell\mid n,\,\ell'\mid n\text{ and }\ell\neq\ell'}\bigl|\Res(\Phi_{F,\ell}^{**},\Phi_{F,\ell'}^{**})\bigr|\\
=\prod_{(\ell,\ell')\in\bN^2:\,\ell\mid\ell'\mid n\text{ and }\ell<\ell'}\bigl|\Res\bigl(\Phi_{F,\ell}^{**},\Phi_{F,\ell'}^{**}\bigr)\bigr|^2
\times\prod_{(\ell,\ell')\in\bN^2:\,\ell\mid n,\,\ell'\mid n\text{ and }\ell\nmid\ell'\nmid\ell}\bigl|\Res\bigl(\Phi_{F,\ell}^{**},\Phi_{F,\ell'}^{**}\bigr)\bigr|.
\end{multline*}
The above three computations yield
\begin{multline*}
\prod_{\ell\,\mid n}\Bigl(\bigl|p_{\ell}^{**}(1;f)\bigr|^\ell\cdot
\prod_{\ell':\,\mid \ell\text{ and }<\ell}\bigl|\Delta_{\ell,\ell'}(f)\bigr|^{\ell'}\Bigr)\\
=\prod_{\ell\,\mid n}\biggl(|\Disc(\Phi_{F,\ell}^{**})|\cdot
\prod_{\ell':\,|\ell\text{ and }<\ell}\bigl|\Res(\Phi_{F,\ell}^{**},\Phi_{F,\ell'}^{**})\bigr|^2
\times\prod_{\ell':\,\ell'\mid n\text{ and }\ell\nmid\ell'\nmid\ell}\bigl|\Res(\Phi_{F,\ell}^{**},\Phi_{F,\ell'}^{**})\bigr| \biggr)
\end{multline*}
(for every $n\in\bN$), and in turn, by M\"obius inversions of both sides, we have
\begin{gather*}
 \bigl|p_n^{**}(1;f)\bigr|^n\cdot\prod_{\ell:\,\mid n\text{ and }<n}\bigl|\Delta_{n,\ell}(f)\bigr|^{\ell}
=|\Disc(\Phi_{F,n}^{**})|\cdot
\prod_{\ell:\,\mid n\text{ and }<n}\bigl|\Res(\Phi_{F,n}^{**},\Phi_{F,\ell}^{**})\bigr|^2,
\end{gather*}
which is, by the definition of $\Delta_{n,n}$, no matter whether $|\Res F|=1$, equivalent to 
\begin{gather}
\bigl|\Delta_{n,n}(f)\bigr|^n
=\biggl|\frac{(\Res F)^{-N_d(n)}\Disc(\Phi_{F,n}^{**})}{\prod_{\ell:\,\mid n\text{ and }<n}(\Delta_{n,\ell}(f))^{n-\ell}}\biggr|
\cdot
\prod_{\ell:\,\mid n\text{ and }<n}
\biggl|\frac{(\Res F)^{-N_d(n,\ell)}\Res(\Phi_{F,n}^{**},\Phi_{F,\ell}^{**})}{\bigl(\Delta_{n,\ell}(f)\bigr)^{\ell}}\biggr|^2.\label{eq:collidetwo}
\end{gather}

{\bfseries (ii)} By the definition \eqref{eq:homogdiscr} of the homogeneous discriminants, 
the labeling \eqref{eq:labelingformal} of 
the sequence 
$(z_{F,j}^{(n)})_{j=1}^{d_n}$, and the
definition \eqref{eq:cyclemutual} of $\Lambda_{F,r}^{(n)}$,
we have
\begin{align*}
 \Disc(\Phi_{F,n}^{**})
=&\prod_{r=0}^{\frac{d_n}{n}-1}\prod_{k=0}^{n-1}\biggl(\prod_{\ell\in\{0,\ldots,n-1\}\setminus\{k\}}(Z_{F,1+nr+k}^{(n)}\wedge Z_{F,1+nr+\ell}^{(n)})\times\prod_{s\neq r}\prod_{\ell=0}^{n-1}(Z_{F,1+nr+k}^{(n)}\wedge Z_{F,1+ns+\ell}^{(n)})\biggr)\\
=&\prod_{r=0}^{\frac{d_n}{n}-1}\Disc\bigl(\Lambda_{F,r}^{(n)}\bigr)\times
\prod_{r=0}^{\frac{d_n}{n}-1}\prod_{s\neq r}\prod_{k=0}^{n-1}\prod_{\ell=0}^{n-1}(Z_{F,1+nr+k}^{(n)}\wedge Z_{F,1+ns+\ell}^{(n)})\\
=&\prod_{r=0}^{\frac{d_n}{n}-1}\Disc\bigl(\Lambda_{F,r}^{(n)}\bigr)\times
\prod_{r=0}^{\frac{d_n}{n}-1}\prod_{s\neq r}\Res\bigl(\Lambda_{F,r}^{(n)},\Lambda_{F,s}^{(n)}\bigr),
\end{align*}
using \eqref{eq:resultantsubst} repeatedly.
If $|\Res F|=1$, then we also have
\begin{align*}
&\log\biggl|\prod_{r=0}^{\frac{d_n}{n}-1}\Disc\bigl(\Lambda_{F,r}^{(n)}\bigr)\biggr|\\
=&\sum_{r=0}^{\frac{d_n}{n}-1}\sum_{k=0}^{n-1}\Bigr(\sum_{m\in\{0,1,\ldots,n-1\}\setminus\{k\}}\Phi_{g_f}(z_{f,1+nr+k}^{(n)},z_{f,1+nr+m}^{(n)})+n\cdot G^F(Z_{F,1+nr+k}^{(n)})\Bigr)
\quad(\text{by }\eqref{eq:homogkernel})\\
=&\sum_{r=0}^{\frac{d_n}{n}-1}\sum_{k=0}^{n-1}\,\sum_{m=1}^{n-1}\Phi_{g_f}\bigl(z_{f,1+nr+k}^{(n)},f^m(z_{f,1+nr+k}^{(n)})\bigr)+0
\quad(\text{by }\eqref{eq:Greenexact}\text{ and the labeling \eqref{eq:labelingformal} of }(z_{F,j}^{(n)}))\\
=&\sum_{r=0}^{\frac{d_n}{n}-1}\sum_{k=0}^{n-1}\,\sum_{m=1}^{n-1}\Bigl(\log|Z_{F,1+nr+k}^{(n)}\wedge F^m(Z_{F,1+nr+k}^{(n)})|-G^F\bigl(Z_{F,1+nr+k}^{(n)}\bigr)-G^F\bigl(F^m(Z_{F,1+nr+k}^{(n)})\bigr)\Bigr)\\
=&\sum_{r=0}^{\frac{d_n}{n}-1}\sum_{k=0}^{n-1}\,\sum_{m=1}^{n-1}\,\sum_{\ell\,\mid m}\log|\Phi_{F,\ell}^{**}(Z_{F,1+nr+k}^{(n)})|
-\Bigl(\sum_{m=1}^{n-1}(d^m+1)\Bigr)\sum_{j=1}^{d_n}G^F\bigl(Z_{F,j}^{(n)}\bigr)\\
=&\log\biggl|\prod_{m=1}^{n-1}\,\prod_{\ell\,\mid m}\Res(\Phi_{F,n}^{**},\Phi_{F,\ell}^{**})\biggr|-0
\quad(\text{by }\eqref{eq:resultantsubst}\text{ repeatedly and }\eqref{eq:Greenexact})\\
=&\log\prod_{\ell:\,\mid n\text{ and }<n}\bigl|\Res(\Phi_{F,n}^{**},\Phi_{F,\ell}^{**})\bigr|^{\frac{n}{\ell}-1}
+\log\prod_{\ell:\,\nmid n\text{ and }<n}\bigl|\Res(\Phi_{F,n}^{**},\Phi_{F,\ell}^{**})\bigr|^{\lfloor\frac{n}{\ell}\rfloor},
\end{align*}
where the third equality is by \eqref{eq:homogkernel}, and
the fourth one is by 
$F^m(Z)\wedge Z=\prod_{\ell\,\mid m}\Phi_{F,\ell}^{**}(Z)$
and \eqref{eq:pullbackGreen}. Hence no matter whether $|\Res F|=1$, we have
\begin{multline}
\biggl|\frac{(\Res F)^{-N_d(n)}\Disc(\Phi_{F,n}^{**})}{\prod_{\ell:\,\mid n\text{ and }<n}(\Delta_{n,\ell}(f))^{n-\ell}}\biggr|\\
=\prod_{\ell:\,\mid n\text{ and }<n}\biggl|\frac{(\Res F)^{-N_d(n,\ell)}\Res(\Phi_{F,n}^{**},\Phi_{F,\ell}^{**})}{\bigl(\Delta_{n,\ell}(f)\bigr)^{\ell}}\biggr|^{\frac{n}{\ell}-1}
\times\prod_{\ell:\,\nmid n\text{ and }<n}\bigl|(\Res F)^{-N_d(n,\ell)}\Res(\Phi_{F,n}^{**},\Phi_{F,\ell}^{**})\bigr|^{\lfloor\frac{n}{\ell}\rfloor}\\
\times\prod_{r=0}^{\frac{d_n}{n}-1}\prod_{s\neq r}
\bigl|(\Res F)^{-M_d(n)}\Res(\Lambda_{F,r}^{(n)},\Lambda_{F,s}^{(n)})\bigr|.\label{eq:collidetwodisc}
\end{multline}
{\bfseries (iii)} Now the desired \eqref{eq:factoredmult} holds by 
\eqref{eq:collideprecise}, 
\eqref{eq:collidetwo}, \eqref{eq:collidetwodisc},
and
$R_{n,\ell}^{**}\in(\bZ[\Rat_d])^*(=\bZ^*=\{\pm 1\})$
if $\ell\nmid n\text{ and }\ell<n$.

\subsection{}\label{sec:conclusion}
The two equalities in Theorem \ref{th:exact} have already been shown in 
Subsections 4.1 and 4.2.

Finally, 
pick
$\tilde{F}=\tilde{F}(X,Y;f)=:\tilde{F}_f(X,Y)
\in((\cO_{\Rat_d(\bC),\an}(D))[X,Y]_d)^2$
on any simply connected domain $D$ in 
the complex manifold $\Rat_d(\bC)$ such that
for any individual $f\in D$, 
$\tilde{F}_f$ is a lift (to $\bC^2$)
of $f$. Then there is a (finitely sheeted) possibly branched 
analytic covering $\eta:M\to D$ of $D$ by 
a complex manifold $M$ such that for every $j\in\{1,\ldots,d_n\}$, 
the mapping $Z_{\tilde{F}_\eta,j}^{(n)}:M\to\bC^2\setminus\{(0,0)\}$ is
complex analytic(, or informally,
$Z_{\tilde{F}_\eta,j}^{(n)}$ is a marked point in $\bC^2\setminus\{(0,0)\}$ 
complex analytically parametrized by $M$). Hence 
recalling the definition \eqref{eq:cyclemutual} of
$\Lambda_{\tilde{F}_\eta,r}^{(n)}$,
for any distinct $r,s\in\{1,\ldots,\frac{d_n}{n}-1\}$, 
$\Res\bigl(\Lambda_{\tilde{F}_\eta,r}^{(n)},\Lambda_{\tilde{F}_\eta,s}^{(n)}\bigr)$
is a (complex analytic so) locally bounded function on $M$, and then
the second equality in \eqref{eq:factoredmult} yields
$\Delta_{n,n}\in\bZ[\Rat_d]$. Now the proof of Theorem \ref{th:exact}
is complete. \qed

\begin{remark}\label{th:integer}
It also follows from an argument similar to that in the final paragraph
in Subsection \ref{sec:conclusion} that $\Disc(\Phi_n^{**})$ is in $\bZ[a,b]$
(rather than merely in the field $\bZ(a,b)$). 
\end{remark}

\begin{remark}
All the potential theoretic computations from Subsection \ref{sec:complexpotential}
and used in the proof of Theorem \ref{th:exact} 
still work over an algebraically closed field $K$ that is complete with respect to
a non-trivial and possibly non-archimedean absolute value $|\cdot|$ (e.g., over
the field $\bC_p$ equipped with the (extended) $p$-adic norm $|\cdot|_p$, or 
the field $\bL$ of ``formal'' Puiseux series $\phi$
centered at the origin $z=0$ in $\bC$
equipped with the norm induced by the orders of $\phi$ at $z=0$).
\end{remark}

\section{A few computations in $\Rat_2$}
\label{sec:quadrat}

Focusing on the case where $d=2$,
we conclude with a few computation in $\Rat_2$.
Here we write each quadratic rational function $f$ on $\bP^1$ as 
$[aX^2+bXY+cY^2:pX^2+qXY+rY^2]=[a:b:c:p:q:r]\in\Rat_2$. 

\subsection{the $n=1$ case} 
When $n=1$, the hypersurface $V(\Delta_{1,1})$
coincides with the hypersurface $V(\Disc(\Phi_1^{**}))$ in 
$\Rat_2=\bP^5\setminus V(\rho_2)$; they
are the loci in $\Rat_2$ for every $[a:b:c:p:q:r]$ 
in which, at least one of the $2_1=2+1=3$ fixed points in $\bP^1$ 
of $f=f_{[a:b:c:p:q:r]}$ is multiple (when working in an algebraically closed
field). We indeed compute
\begin{gather*}
 \rho_2=a^{2} r^{2} - a b q r - 2 a c p r + a c q^{2} + b^{2} p r - b c p q + c^{2} p^{2}
\end{gather*}
and
\begin{dmath*}
 \Disc(\Phi_1^{**})=4 a^{3} c - a^{2} b^{2} + 2 a^{2} b r - 12 a^{2} c q - a^{2} r^{2} + 2 a b^{2} q + 18 a b c p - 4 a b q r - 18 a c p r + 12 a c q^{2} + 2 a q r^{2} - 4 b^{3} p + 12 b^{2} p r - b^{2} q^{2} - 18 b c p q - 12 b p r^{2} + 2 b q^{2} r + 27 c^{2} p^{2} + 18 c p q r - 4 c q^{3} + 4 p r^{3} - q^{2} r^{2},
\end{dmath*}
so that, by the first equality in Theorem \ref{th:exact},
\begin{dmath*}
 \Delta_{1,1}=\pm D_1^{**}=\pm(\rho_2)^{-1}\Disc(\Phi_1^{**})
=\pm(4 a^{3} c - a^{2} b^{2} + 2 a^{2} b r - 12 a^{2} c q - a^{2} r^{2} + 2 a b^{2} q + 18 a b c p - 4 a b q r - 18 a c p r + 12 a c q^{2} + 2 a q r^{2} - 4 b^{3} p + 12 b^{2} p r - b^{2} q^{2} - 18 b c p q - 12 b p r^{2} + 2 b q^{2} r + 27 c^{2} p^{2} + 18 c p q r - 4 c q^{3} + 4 p r^{3} - q^{2} r^{2})/\rho_2
\end{dmath*}
on $\Rat_2$.

\subsection{the $n=2$ case} When $n=2$, since $2_2=(2^2+1)-(2^1+1)=2$,
for each $[a:b:c:p:q:r]\in\Rat_2$,
there are exactly one cycle $C_2$ in $\bP^1$ of $f=f_{[a:b:c:p:q:r]}$ 
having the formally exact periods $n=2$, and the hypersurface $V(\Disc(\Phi^{**}_2))$
coincides with the hypersurface $V(\Delta_{2,1})$ in $\Rat_2$; they are
the loci in $\Rat_2$ for every $[a:b:c:p:q:r]$ in which, 
this cycle $C_2$ in $\bP^1$ of $f=f_{[a:b:c:p:q:r]}$ 
reduces to a fixed point of $f$ (when working in an algebraically closed field).
We indeed compute
\begin{dmath*}
 \Disc(\Phi_2^{**})=4 a^{3} c - a^{2} b^{2} + 2 a^{2} b r + 4 a^{2} c q + 3 a^{2} r^{2} - 2 a b^{2} q + 2 a b c p + 6 a c p r + 2 a q r^{2} - b^{2} q^{2} + 2 b c p q + 4 b p r^{2} - 2 b q^{2} r - c^{2} p^{2} + 2 c p q r + 4 p r^{3} - q^{2} r^{2}
\end{dmath*}
and
\begin{dmath*}
 \Res(\Phi_2^{**},\Phi_1^{**})=
4 a^{5} c r^{2} - a^{4} b^{2} r^{2} - 4 a^{4} b c q r + 2 a^{4} b r^{3} - 8 a^{4} c^{2} p r + 4 a^{4} c^{2} q^{2} + 4 a^{4} c q r^{2} + 3 a^{4} r^{4} + a^{3} b^{3} q r + 6 a^{3} b^{2} c p r - a^{3} b^{2} c q^{2} - 4 a^{3} b^{2} q r^{2} - 4 a^{3} b c^{2} p q - 2 a^{3} b c p r^{2} - 2 a^{3} b c q^{2} r - 3 a^{3} b q r^{3} + 4 a^{3} c^{3} p^{2} - 8 a^{3} c^{2} p q r + 4 a^{3} c^{2} q^{3} + 3 a^{3} c q^{2} r^{2} + 2 a^{3} q r^{4} - a^{2} b^{4} p r + a^{2} b^{3} c p q + 2 a^{2} b^{3} p r^{2} + 2 a^{2} b^{3} q^{2} r - a^{2} b^{2} c^{2} p^{2} + 4 a^{2} b^{2} c p q r - 2 a^{2} b^{2} c q^{3} + 3 a^{2} b^{2} p r^{3} - a^{2} b^{2} q^{2} r^{2} - 2 a^{2} b c^{2} p^{2} r - 2 a^{2} b c^{2} p q^{2} - 7 a^{2} b c p q r^{2} + 4 a^{2} b p r^{4} - 4 a^{2} b q^{2} r^{3} + 4 a^{2} c^{3} p^{2} q - 10 a^{2} c^{2} p^{2} r^{2} + 6 a^{2} c^{2} p q^{2} r - 2 a^{2} c p q r^{3} + 2 a^{2} c q^{3} r^{2} + 4 a^{2} p r^{5} - a^{2} q^{2} r^{4} - 2 a b^{4} p q r + 2 a b^{3} c p^{2} r + 2 a b^{3} c p q^{2} + a b^{3} q^{3} r - 4 a b^{2} c^{2} p^{2} q + 6 a b^{2} c p^{2} r^{2} - a b^{2} c q^{4} - 2 a b^{2} p q r^{3} + 2 a b^{2} q^{3} r^{2} + 2 a b c^{3} p^{3} - 9 a b c^{2} p^{2} q r + 2 a b c^{2} p q^{3} - 8 a b c p^{2} r^{3} + 4 a b c p q^{2} r^{2} - 2 a b c q^{4} r - 4 a b p q r^{4} + a b q^{3} r^{3} + 8 a c^{3} p^{3} r - a c^{3} p^{2} q^{2} - 2 a c^{2} p^{2} q r^{2} + 2 a c^{2} p q^{3} r - 8 a c p^{2} r^{4} + 6 a c p q^{2} r^{3} - a c q^{4} r^{2} - b^{4} p q^{2} r + 2 b^{3} c p^{2} q r + b^{3} c p q^{3} + 4 b^{3} p^{2} r^{3} - 2 b^{3} p q^{2} r^{2} - b^{2} c^{2} p^{3} r - 3 b^{2} c^{2} p^{2} q^{2} - 2 b^{2} c p^{2} q r^{2} + 2 b^{2} c p q^{3} r + 4 b^{2} p^{2} r^{4} - b^{2} p q^{2} r^{3} + 3 b c^{3} p^{3} q + 4 b c^{2} p^{3} r^{2} - 4 b c^{2} p^{2} q^{2} r - 4 b c p^{2} q r^{3} + b c p q^{3} r^{2} - c^{4} p^{4} + 2 c^{3} p^{3} q r + 4 c^{2} p^{3} r^{3} - c^{2} p^{2} q^{2} r^{2}
=\rho_2\cdot\Disc(\Phi_2^{**}),
\end{dmath*}
so that, also by the equalities in Theorem \ref{th:exact}, 
\begin{dmath*}
\Delta_{2,1}=\pm R_{2,1}^{**}
=\pm(\rho_2)^{-2}\Res(\Phi_2^{**},\Phi_1^{**})
=\pm(\rho_2)^{-1}\Disc(\Phi_2^{**})
=\pm D_2^{**}\\
=\pm(4 a^{3} c - a^{2} b^{2} + 2 a^{2} b r + 4 a^{2} c q + 3 a^{2} r^{2} - 2 a b^{2} q + 2 a b c p + 6 a c p r + 2 a q r^{2} - b^{2} q^{2} + 2 b c p q + 4 b p r^{2} - 2 b q^{2} r - c^{2} p^{2} + 2 c p q r + 4 p r^{3} - q^{2} r^{2})/\rho_2
\end{dmath*}
on $\Rat_2$ and
\begin{gather*}
(\Delta_{2,2})^2=\pm\frac{D_2^{**}}{\Delta_{2,1}}=1\quad\text{or equivalently}\quad\Delta_{2,2}=\pm 1
\end{gather*}
on $\Rat_2$.

\subsection{the $n=3$ case} When $n=3$, we similarly compute 
\begin{dmath*}
\Delta_{3,1}=\pm R_{3,1}^{**}
=\pm(\rho_2)^{-6}\Res(\Phi^{**}_3,\Phi^{**}_1)\quad(\text{by the second equality in Theorem \ref{th:exact}})
=\pm(16 a^{6} c^{2} - 8 a^{5} b^{2} c + 16 a^{5} b c r + 4 a^{5} c r^{2} + a^{4} b^{4} - 4 a^{4} b^{3} r - 8 a^{4} b^{2} c q + 3 a^{4} b^{2} r^{2} + 48 a^{4} b c^{2} p + 4 a^{4} b c q r + 2 a^{4} b r^{3} + 24 a^{4} c^{2} p r + 12 a^{4} c^{2} q^{2} + 16 a^{4} c q r^{2} + 7 a^{4} r^{4} + 2 a^{3} b^{4} q - 20 a^{3} b^{3} c p - 5 a^{3} b^{3} q r + 30 a^{3} b^{2} c p r - 11 a^{3} b^{2} c q^{2} - 3 a^{3} b^{2} q r^{2} - 12 a^{3} b c^{2} p q + 18 a^{3} b c p r^{2} - 2 a^{3} b c q^{2} r - 5 a^{3} b q r^{3} + 36 a^{3} c^{3} p^{2} + 48 a^{3} c^{2} p q r - 8 a^{3} c^{2} q^{3} + 8 a^{3} c p r^{3} + 7 a^{3} c q^{2} r^{2} + 2 a^{3} q r^{4} + 2 a^{2} b^{5} p - 7 a^{2} b^{4} p r + 3 a^{2} b^{4} q^{2} - 9 a^{2} b^{3} c p q + 5 a^{2} b^{3} p r^{2} - 3 a^{2} b^{3} q^{2} r + 27 a^{2} b^{2} c^{2} p^{2} - 18 a^{2} b^{2} c p q r - a^{2} b^{2} c q^{3} + 7 a^{2} b^{2} p r^{3} + 54 a^{2} b c^{2} p^{2} r + 18 a^{2} b c^{2} p q^{2} + 63 a^{2} b c p q r^{2} - 16 a^{2} b c q^{3} r + 16 a^{2} b p r^{4} - 3 a^{2} b q^{2} r^{3} + 54 a^{2} c^{2} p^{2} r^{2} - 18 a^{2} c^{2} p q^{2} r + 9 a^{2} c^{2} q^{4} + 18 a^{2} c p q r^{3} + 5 a^{2} c q^{3} r^{2} + 4 a^{2} p r^{5} + 3 a^{2} q^{2} r^{4} + 2 a b^{5} p q - 12 a b^{4} c p^{2} - 4 a b^{4} p q r + 2 a b^{4} q^{3} + 12 a b^{3} c p^{2} r - 12 a b^{3} c p q^{2} - 16 a b^{3} p q r^{2} + a b^{3} q^{3} r - 27 a b^{2} c^{2} p^{2} q - 18 a b^{2} c p^{2} r^{2} - a b^{2} c q^{4} - 2 a b^{2} p q r^{3} - 3 a b^{2} q^{3} r^{2} + 54 a b c^{3} p^{3} + 81 a b c^{2} p^{2} q r - 30 a b c^{2} p q^{3} + 48 a b c p^{2} r^{3} - 18 a b c p q^{2} r^{2} - 4 a b c q^{4} r + 4 a b p q r^{4} - 5 a b q^{3} r^{3} + 27 a c^{3} p^{2} q^{2} + 54 a c^{2} p^{2} q r^{2} + 12 a c^{2} p q^{3} r - 3 a c^{2} q^{5} + 24 a c p^{2} r^{4} + 30 a c p q^{2} r^{3} - 7 a c q^{4} r^{2} + 16 a p q r^{5} - 4 a q^{3} r^{4} + b^{6} p^{2} - 3 b^{5} p^{2} r + 2 b^{5} p q^{2} + 3 b^{4} c p^{2} q + 9 b^{4} p^{2} r^{2} - b^{4} p q^{2} r + b^{4} q^{4} - 9 b^{3} c^{2} p^{3} - 30 b^{3} c p^{2} q r + 5 b^{3} c p q^{3} - 8 b^{3} p^{2} r^{3} - b^{3} p q^{2} r^{2} + 2 b^{3} q^{4} r + 27 b^{2} c^{2} p^{3} r + 18 b^{2} c p^{2} q r^{2} - 12 b^{2} c p q^{3} r + 2 b^{2} c q^{5} + 12 b^{2} p^{2} r^{4} - 11 b^{2} p q^{2} r^{3} + 3 b^{2} q^{4} r^{2} - 27 b c^{3} p^{3} q - 27 b c^{2} p^{2} q^{2} r + 3 b c^{2} p q^{4} - 12 b c p^{2} q r^{3} - 9 b c p q^{3} r^{2} + 2 b c q^{5} r - 8 b p q^{2} r^{4} + 2 b q^{4} r^{3} + 27 c^{4} p^{4} + 54 c^{3} p^{3} q r - 9 c^{3} p^{2} q^{3} + 36 c^{2} p^{3} r^{3} + 27 c^{2} p^{2} q^{2} r^{2} - 12 c^{2} p q^{4} r + c^{2} q^{6} + 48 c p^{2} q r^{4} - 20 c p q^{3} r^{3} + 2 c q^{5} r^{2} + 16 p^{2} r^{6} - 8 p q^{2} r^{5} + q^{4} r^{4})/(\rho_2)^2
\end{dmath*}
on $\Rat_2$; the computations of $D_3^{**}=(\rho_2)^{-15}\Disc(\Phi^{**}_3)$ and $\Delta_{3,3}$ on $\Rat_2$
are more lengthy.

\begin{acknowledgement}
The author thanks the referee and the editors
for a very careful scrutiny and invaluable comments.
The author was partially supported by JSPS Grant-in-Aid 
for Scientific Research (C), 19K03541 and (B), 19H01798.
\end{acknowledgement}

\def\cprime{$'$}

\end{document}